\documentclass[11pt]{article}
\usepackage{amssymb}
\usepackage{latexsym}
\newtheorem{thm}{Theorem}[section]
\newtheorem{prop}[thm]{Proposition} \newtheorem{lemma}[thm]{Lemma}
\newtheorem{cor}[thm]{Corollary} 
 \newtheorem{rmk}[thm]{Remark}

\newcommand {\pf}{\noindent{\bf Proof.}\ }
\newcommand{\reals}{{\mathbb R}}

\newcommand{\ver}{{\rm ver}}
\newcommand{\hor}{{\rm hor}}
\newcommand{\Hom}{{\rm Hom}}
\newcommand{\Tr}{{\rm Tr}}
\newcommand{\covder}{{\bf D}}
\newcommand{\grad}{\nabla}
\newcommand{\ave}{\cala}
\newcommand{\doub}{\calv}
\newcommand{\gammabar}{{\overline{\Gamma}}}
\newcommand{\fbar}{{\overline{F}}}
\newcommand{\nbar}{{\overline{N}}}
\newcommand{\pbar}{{\overline{P}}}
\newcommand{\proj}{{\mathbb P}}
\newcommand{\cala}{{\cal A}}
\newcommand{\cald}{{\cal D}}
\newcommand{\calh}{{\cal H}}

\newcommand{\calt}{{\cal T}}
\newcommand{\calv}{{\cal V}}

\newcommand{\qed}{\begin{flushright} $\Box$\ \ \ \ \ \end{flushright}}
\newcommand{\lqed}{\begin{flushright} $\triangle$\ \ \ \ \ \end{flushright}}

\title{{\bf Almost Invariant Submanifolds for
Compact Group Actions }} \author{Alan
Weinstein\thanks{Research partially supported by NSF Grants
DMS-96-25122 and DMS-99-71505
and the Miller Institute for Basic Research in Science.}
\\Department of Mathematics\\ University of California\\ Berkeley, CA
94720 USA\\ {\small(alanw@math.berkeley.edu)}}
\date{July 29, 1999}
\begin{document}
\setlength{\baselineskip}{15pt}

\maketitle

\section{Introduction}
\label{sec-intro}
A compact (not necessarily connected) Lie group $G$
carries a (unique) bi-invariant
probability measure.  Using this measure, one can average orbits of
actions of $G$ on affine convex sets to obtain fixed points.  In
particular, if $G$
acts on a manifold $M$, $G$ leaves invariant a riemannian metric on
$M$, and this metric can sometimes be used to obtain fixed points for
the nonlinear action of $G$ on $M$ itself.  By this method, \'Elie Cartan
proved that a compact group $G$ acting by
isometries on a simply-connected manifold $M$ of nonpositive sectional
curvature always has a fixed point.\footnote{In Note III, $\rm{n^o}$ 19
  of \cite{ca:lecons}, Cartan proved the fixed point theorem for a
  finite group; he then remarked on p. 19 of \cite{ca:groupes} that the
  same argument works for a compact group.  See Theorem 13.5 in
  Chapter 1 of Helgason's book \cite{he:differential} for a full proof
  in the compact case.}
Cartan's result was extended by Grove and
Karcher \cite{gr-ka:how} to arbitrary
manifolds under the assumption that the
$G$-action has an orbit which is sufficiently small relative
to a distance scale provided by the geometry of $M$.

In this paper, we develop a method for averaging nearby submanifolds
in a riemannian manifold.  This enables us to extend the Grove-Karcher
theorem from points to submanifolds; i.e. we establish that, if a riemannian
$G$-manifold has a compact submanifold $N$ whose images under the
$G$-action are sufficiently $C^1$-close to one another, relative to the
geometry of the pair $(M,N)$, then there is a $G$-invariant
submanifold near $N$. Extending the Cartan-Grove-Karcher nonlinear
averaging method to submanifolds requires new estimates on
the geometry of tubular neighborhoods.  Our construction is based in
an essential way on Whitney's
idea\footnote{See Chapter V of the cited article.  After finishing
this paper, I looked again at this reference and discovered that
Whitney used a center of mass construction, too!}
\cite{wh:differentiable} of realizing submanifolds as zeros of
sections of vector bundles.

Our result fits nicely into the ``stability'' framework promoted by
Ulam (see for example page 584 of \cite{ul:sets}, as well as Anderson's
interesting discussion using nonstandard analysis in
\cite{an:almost}).  Ulam asks when a mathematical object which
``almost solves'' a certain problem can be shown to be near an
exact solution.  In fact, our interest in
the invariant submanifold problem was
motivated by some cases of the following ``almost-homomorphism
problem'': given a map $\phi:G\to H$ between groups such that
$\phi(gh)$ is close in some sense to $\phi(g)\phi(h)$ for all $g$
and $h$ in $G$, can one conclude that there is a homomorphism near
$\phi$?  This problem was posed by Ulam and solved by Hyers
\cite{hy:stability} for $G=H=\reals$. It was solved for compact Lie
groups $G$ and $H$, with effective estimates, by Grove,
Karcher, and Ruh \cite{gr-ka-ru:group} using the Grove-Karcher
fixed-point theorem.   When $G$ is compact and $H$ is the unitary
group of a Hilbert space, it was solved by
de la Harpe and Karoubi \cite{de-ka:representations}.

Our invariant submanifold theorem implies an
almost-homomorphism theorem when $G$ is finite and $H$ is the
group of diffeomorphisms of a compact manifold.
Our ultimate aim
is to extend the almost-homomorphism theorem to
cases where $G$ is a
compact group and $H$ is the group of smooth bisections of a Lie groupoid.
(See \cite{ca-we:lectures}.)  Such
an extension would imply very strong structure
theorems for certain groupoids.  We refer to the lecture \cite{we:from} for
an overview
of these potential applications and their connection with symplectic
geometry.   Details of these applications will appear in subsequent papers.

An important step in our proof requires averaging in a Grassmann
manifold. This can be accomplished by Cartan-Grove-Karcher
averaging, but we get estimates better suited to our purposes by
identifying subspaces with projections and averaging the
projections. The averaging of projections has already been
treated, in finite and infinite dimensions, by a number of authors
over several decades, including de la Harpe and Karoubi in the paper
cited above.
We review some of this work and derive the estimates we need in an Appendix.

The body of the paper is organized as follows.  In Section 2, we define
a $C^1$-distance between submanifolds of a riemannian manifold and, in
terms of
this distance, state our main theorem on invariant manifolds, Theorem
\ref{thm-main}, which
follows immediately from the averaging theorem, Theorem \ref{thm-centerofmass}.
Section 3 develops estimates on geometric objects in tubular
neighborhoods; the reader may wish to skim this section and return
to it after reading Section 4, where we first outline the proof
of the averaging theorem and then fill in the details.
In Section 5, we discuss some alternate approaches to the
averaging of submanifolds which may be interesting in their own
right, though we have not yet succeeded in implementing any of
them.  We conclude in Section 6 with a discussion
possible extensions of our results, and the relation of our work to
other areas.

In view of potential applications to almost actions of compact Lie
groups, we note that most of
our work extends without modification to the case where $M$ is a Hilbert
manifold, though the
submanifold $N$ must remain finite-dimensional.   There is just one
issue, the application of degree theory, which now restricts our main
result to finite-dimensional $M$.

We have not tried to find optimal values for the many constants in the estimates
in Sections \ref{sec-geometry} and \ref{sec-main}.  Rather than
denoting these constants by letters, though, we have chosen
numerical values which are valid and which may be recognized when they
reappear.  The resulting constants in
the main theorems are therefore very far from optimal.
The reader is invited to find better ones.

A word on notation: we will write $d(\cdot,\cdot)$
to denote all of the many distance functions used in this paper, decorating
it with subscripts only when necessary to avoid ambiguity.  We will
use bars $|\cdot|$ for lengths of paths and vectors, and double
bars $\|\cdot\|$ for operator norms.

Many people have offered useful suggestions in the course of this
work.  I would particularly like to thank David Aldous, Eugenio Calabi,
Marc Chamberland,
Jost Eschenburg, Robert Greene, Mikhail Gromov, Karsten Grove, Hermann Karcher,
Yi Ma, Yoshiaki Maeda, and Jamie Sethian.

\section{Bounded geometry and distance between submanifolds}
Let $N$ be a closed submanifold of the riemannian manifold $M$.  The
{\bf normal bundle} $\nu N \subset TM$ consists of those tangent
vectors along $N$ which are orthogonal to $TN$.  The open ball bundle
of radius $r$ in $\nu N$ will be denoted by $\nu^r N$.

We define the {\bf normal exponential map} $\exp_N$ to be
the restriction to $\nu N$ of the exponential map $\exp:TM\to M.$
(If $M$ is not complete, $\exp_N$, like $\exp$, is defined only on
a proper open neighborhood of the zero section.)
The {\bf normal
injectivity radius} $i_N$ of $N$ is
$$\sup \{r|\exp_N ~\mbox{(is defined and) is an embedding on}~ \nu^r
N\}.$$ (When $N$ is a point, its normal injectivity radius is just the
usual injectivity radius of $M$ at that point.) The restriction of
$\exp_N$ to $\nu^{i_N} N$ is still an embedding; we call its image
$\calt N$ {\em the} tubular neighborhood of $N$.  For any $r \in
(0,i_N]$, the set
$$\calt^r N =\{x\in M|d(x,M)<r\}=\exp(\nu^r N)$$ will be called the $r$-{\bf
tube} about $N$.

The bundle retraction in $\nu N$ is transferred by $\exp_N$ to a
retraction $\pi_N:\calt N\rightarrow N$ which maps each point of
$\calt N$ to the nearest point in $N$.  The fibres of $\pi_N$ will be
called the {\bf normal slices} for $N$.

\subsection{The Grassmann bundle}
We will denote by $GM$ the Grassmann bundle over $M$, whose fibre over each
point
$x$ is the (disconnected) manifold of all finite-codimensional subspaces of
the tangent space $T_xM$.
 We give each
fibre the Finsler metric described in the
Appendix.

 These fibre metrics define a supremum distance
on the space of sections of $GM$ over any subset of $M$.

The {\bf Gauss map} of a submanifold $N$ is defined as the section of
$GM$ over $N$ whose value at each point is the   normal  space to
$N$.  This section extends in a natural way to a section $\Gamma_
N:\calt N\rightarrow GM$ defined over the tubular neighborhood of $N$:
we take the normal spaces to $N$ and parallel translate them along all
the geodesics normal to $N$.  The image of $\Gamma _N$, when thought
of as a vector bundle over the tubular neighborhood, will be called
the {\bf vertical bundle} and its
orthogonal complement (obtained by translating the tangent spaces of
$N$) the {\bf horizontal bundle}.

The vertical bundle contains the tangent vectors to the normal geodesics
but does not in
general coincide with the tangent bundle along the normal slices.  We
will refer to the latter as the {\bf quasi-vertical bundle}.
(The two bundles are compared in Section \ref{subsec-twobundles}.)
The orthogonal complement to the quasi-vertical bundle will be called
the {\bf quasi-horizontal bundle}.

\subsection{Bounded geometry and $C^1$ distance}
To define a geometrically meaningful $C^1$ distance
between submanifolds, we should measure distances and angles in the
same units.  Since angles are dimensionless, we
do this by choosing a positive number $c$ which
functions as a ``unit of inverse length.''  Given $c$, we say
that a pair $(M,N)$ consisting of a manifold $M$ and its submanifold
$N$ has {\bf   geometry bounded by} $c$ if: (i) $i_N >
1/c$; (ii) the sectional curvatures of $M$ in the $1/c$-tube about $N$
are bounded in absolute
value by $c^2$; (iii) the injectivity radius of each point in the tube is at
least $1/c$.  If
$N$ is compact, $(M,N)$ necessarily\footnote{Even if $M$ has infinite
dimension!} has   geometry bounded by some $c$.

By rescaling the metric, we can always convert geometry bounded by $c$
to geometry bounded by 1.  (Note that rescaling the metric does not
change the metric on the Grassmann manifold.)  We will therefore
restrict the statement of our results to this case, leaving it to the
reader who so desires to put the $c$'s back in.  A pair $(M,N)$ with
geometry bounded by $1$ will be called a {\bf gentle} pair.

\begin{rmk}
\label{rmk-condition}
{\em Condition (iii) in the definition of bounded
geometry is used in only a couple of places below.  Perhaps
 our results can be obtained without
this condition.  One example where the condition is violated is
a geodesic circle $N$ in a very thin cylinder
$M$.}
\end{rmk}

If $N$ and $N'$ are submanifolds of $M$ such that $N'$ lies in
$\calt N$ and is the image under $\exp_N$ of a section
$\sigma(N,N')$ of $\nu N$, we define their (nonsymmetric! but see
Remark \ref{rmk-quasisymmetry}) $C^1$-distance $d(N,N')$ as follows.   We
assign two numbers to each to each $x'\in N'$: the
length of the geodesic segment $\tau$ from $x'$ to the nearest
point $x$ on $N$ and the distance between $T_{x'}N'$ and the
parallel translate of $T_x N$ along $\tau$. The $C^1$-distance
$d(N,N')$ is defined as the supremum of all
these numbers as $x'$ ranges over $N'$.  Since
the distance between
two subspaces equals the distance between their orthogonal
complements (Corollary \ref{cor-duality}), we can use normal spaces
instead of tangent spaces and interpret $d(N,N')$ as
the maximum of the supremum
norm of $\sigma(N,N')$ and the distance between the extended Gauss
map $\Gamma_N$ and the Gauss map $\Gamma_{N'}$ as sections of $GM$
over $N'$.

We will occasionally use a $C^0$ distance as well, defined simply of
the maximum over $N'$ of the distance to the nearest point of $N$.
This distance is defined for any pair of compact subsets of a metric
space.

We can now state the main result of this paper.

\begin{thm}
\label{thm-main}
Let $M$ be a riemannian $G$-manifold for a compact group $G$, and
let
$(M,N)$ be a gentle pair.
If $d(N,gN) < \epsilon < \frac{1}{20000}$ for all $g \in G$, then there is
a $G$-invariant
submanifold $\overline N$ with $d(N,\overline {N}) < 136 \sqrt{\epsilon}$.
\end{thm}

This theorem follows immediately from a more general statement,
analogous to Grove and Karcher's center of mass construction
\cite{gr-ka:how} for points.

\begin{thm}
\label{thm-centerofmass}
Let $M$ be a riemannian manifold and $\{N_g\}$ a family of compact
submanifolds of $M$ parametrized in a measurable way\footnote{
We will not attempt to give a precise definition of this notion,
which should involve measurable dependence of derivatives on
parameters to permit the exchange of integration and differentiation.
In the applications we have in mind, the parameter
space is a manifold and the dependence is smooth.
} by elements $g$ of a probability space $G$, such that all
the pairs $(M,N_g)$ are gentle.  If $d(N_g,N_h) <
\epsilon < \frac{1}{20000}$ for all $g$ and $h$ in $G$, there is a well
defined {\bf
center of mass} submanifold $\overline N$ with $d(N_g,\overline {N}) <
136\sqrt{\epsilon}$ for all $g$ in $G$.
The center of mass construction is equivariant with respect to
isometries of $M$ and measure preserving automorphisms of $G$.
\end{thm}

\begin{rmk}
{\em When $G$ consists of two points with equal masses, we obtain a natural
midpoint construction for pairs of nearby submanifolds in a riemannian
manifold.  We know of no simpler construction, even for pairs of
embedded circles in euclidean space.}
\end{rmk}

The numbers $20000$ and $136$ in the theorems above can certainly be
reduced if more care is taken in the many estimates used in the
proof.   On the other hand, it is not at clear whether the factor
$\sqrt{\epsilon}$ can be replaced by a multiple of $\epsilon$ itself.
Our proof does show that the $C^0$ distance from $N$ (or $N_g$) to
$\nbar$ is bounded by $100\epsilon$, but we do not see how to get an
estimate linear in $\epsilon$ for the derivative.   (See Remark
\ref{rmk-squareroot}.)

\section{Geometry in  tubular neighborhoods}
\label{sec-geometry}

In this section, we will establish some
geometric estimates needed for proving the main theorem.

\subsection{Estimate on the second fundamental form}

Recall that the second fundamental
form $B$ of $N$ is defined for a normal vector $v$ and tangent vectors
$w_i$ at $p\in N$ by
\[B_v(w_1,w_2)=-\langle v,\covder_{w_1}W_2\rangle ,\]
where $W_2$ is a vector field tangent to $N$ with $W_2(p)=w_2$,
$\covder$ is the riemannian connection, and $\langle~,~\rangle$ is the
riemannian inner product.\footnote{We follow the sign convention used
by Warner \cite{wa:extensions} (following \cite{bi-cr:geometry}) and
Eschenburg \cite{es:comparison}.  This convention is consistent with
the interpretation of the second fundamental form as the ``first
variation of the metric,'' but most authors use the opposite sign,
consistent with the direction of principal curvature vectors of curves
in a submanifold.  Lang \cite{la:fundamentals} introduces
second fundamental forms with {\em both} signs as the off-diagonal
entries in a single $2\times 2$ block matrix.}  The norm
$\|B\|$ is then defined as the supremum of $|B_v(w,w)|$ as $v$ and $w$
range over unit normal and tangent vectors respectively, and their
basepoint ranges over $N$.

Since focal points along geodesics normal
 to a submanifold cannot occur within the normal
injectivity radius, gentleness of $(M,N)$
 implies that all focal points to $N$
occur at a distance greater than $1$ from $N$.  This fact and
Proposition \ref{prop-focal} below will imply a bound for the
second fundamental form of $N$ (Corollary \ref{cor-secondbound}).

\begin{prop}
\label{prop-focal}
Let $v$ and $w$ be unit normal and tangent vectors respectively at $p
\in N$.  Suppose that $B_v(w,w) = -b$ for some $b > 0$ and that the
sectional curvatures of $M$ are bounded above by $-1$.  Then the
geodesic segment $\gamma$ with initial velocity $v$ and length $L$
contains at least one focal point to $N$ as long as $b \geq 1/L + L/2$.
\end{prop}
\pf We use a standard argument to show that the index
form of $\gamma$ (with appropriate boundary conditions) is not
positive definite.

Using parallel translation and a slight abuse of notation, we
represent any vector field $X$ along $\gamma$ by a function
$X:[0,L]\rightarrow T_pM$.  The index form on the space of
such $X$ with $X(t)\in T_pN$ and $X(L)=0$ is defined by the integral:
$$
I(X,X)=B_v(X(0),X(0)) + \int _0^L \left(\langle X'(t),X'(t) \rangle -
\langle R(t)X(t),X(t)\rangle\right)dt,
$$
 where $R(t)$ is the operator given via the Riemann curvature tensor as
$$X\mapsto R_{\gamma(t)}(\gamma'(t), X)\gamma'(t).$$
 For the
``linearly decreasing'' vector field $X(t) = (1-t/L)w$,
\begin{eqnarray*}
I(X,X) & = &B_v(w,w)+\int_0^L((1/L^2) - (1-t/L) \langle
  R(t)w,w\rangle)dt\\ & = &-b + 1/L - \int_0^L (1-t/L)
  \langle R(t)w,w\rangle dt.
\end{eqnarray*}
Since the sectional curvatures are at least $-1$,
$$
I(X,X) \leq -b + 1/L + \int_0^L (1-t/L)dt = -b + 1/L + L/2.
$$
When $b\geq 1/L +  L/2$, $I(X,X) \leq 0$, so there must be at least
one focal point on $\gamma$.  \qed

\begin{cor}
\label{cor-secondbound}
If $(M,N)$ is a gentle pair, then the second
fundamental form of $N$ satisfies the bound $\|B\| \leq \frac{3}{2}$.
\end{cor}
\pf The distance from $N$ to any focal point is at least $1$.  If
$|B_v(w,w)| > \frac{3}{2}$ for some unit vectors $v$ and $w$, we can change
the sign of $v$ if necessary to assure that $B_v(w,w) < -\frac{3}{2}$, and
Proposition \ref{prop-focal} leads to a contradiction.  \qed

\subsection{Comparison of curves on and near $N$}

The following lemma shows that the retraction $\pi_N:\calt
N\rightarrow N$ does not increase lengths too much.
\begin{lemma}
\label{lem-warner}
Let $(M,N)$ be a gentle pair.  Let
 $\sigma$ be a path in the
 $r$-tube about $N$ for some
$r < \frac{1}{2}$.  Then
\begin{equation}
\label{eq-projection}
 (\cos r -
 \textstyle{\frac{3}{2}}\sin r) |\pi_N \circ \sigma| \leq |\sigma|.
\end{equation}
 (The  function in parentheses
decreases with $r$ and takes the values $1$ at $r=0$, $.59\ldots$
 at $r=\frac{1}{4}$, and $.15\ldots$ at $r=\frac{1}{2}$.)
\end{lemma}
\pf Assume that $\sigma$ is defined on the interval $[0,1]$ and
construct a ``rectangle'' $\Sigma:[0,1]\times [0,1] \rightarrow M$ as
follows: for each $s\in [0,1]$, $s\mapsto \Sigma (s,1)$ is the
minimizing geodesic from $\pi_N(\sigma (s))$ to $\sigma(s)$,
parametrized proportionally to arc length.  Corresponding tangent vectors to
$\sigma$ and $\pi_N\circ\sigma$ are thus the terminal and initial
values of $N$-Jacobi fields along normal geodesics.  By Warner's
 generalization
of the Rauch comparison theorem (\cite{wa:extensions}, Theorem 4.3(b)),
using the
upper bound $1$ for the sectional curvature and the bound $B_{\sigma '(0)} \geq
-\frac{3}{2}I$
(Corollary \ref{cor-secondbound}) on the second
fundamental form of $N$, we obtain the inequality
$|\sigma'(s) |\geq (\cos r - \frac{3}{2}\sin r)|(\pi_N\circ\sigma)'(s)|$
for all $s$, which gives (\ref{eq-projection}).
\qed

\begin{cor}
\label{cor-seven}
Let $(M,N)$ be a gentle pair, $x$ and $y$ points of the $\frac{1}{4}$-tube
about $N$ such that
$d(x,y)<\frac{1}{2}$.  Then
$$d(\pi_N(x),\pi_N(y))< 7 d(x,y).$$
\end{cor}
\pf
A path of length less than $\frac{1}{2}$ from $x$ to $y$ remains in the
$\frac{1}{2}$-tube about $N$.
 Lemma \ref{lem-warner} gives the desired inequality  since
$1/.15 < 7$,
\qed

A similar argument gives:

\begin{cor}
\label{cor-mndistance}
For a gentle pair $(M,N)$, let $x$ and $y$ be points of $N$ such that
$d_M(x,y)=\alpha< \frac{1}{2}$.  Then
$$
d_M(x,y) \leq d_N(x,y)\leq \frac{1}{\cos\frac{1}{2}\alpha -
\frac{3}{2}\sin \frac{1}{2} \alpha}d_M(x,y) \leq 7d_M(x,y).
$$
\end{cor}

\begin{rmk}
\label{rmk-surely}
{\em
This estimate and the ones to follow can surely be sharpened, but the
versions which we give suffice for our purposes.
}
\end{rmk}

\subsection{Vertical and quasi-vertical bundles}
\label{subsec-twobundles}

The quasi-vertical bundle, tangent to the normal slices for
$N\subset M$, is generally not parallel along normal geodesics and
thus does not coincide with the vertical bundle. We will estimate
the distance between the two bundles by estimating the covariant
derivative of the quasi-vertical bundle along normal geodesics.

\begin{prop}
\label{prop-twobundles}
For $x\in \calt N$, the distance between the quasi-vertical and
vertical spaces is at most $\frac{1}{4} d(x,N)^2$.  The distance
between the projections on this spaces is at most $\frac{1}{5}d(x,N)^2.$
\end{prop}
\pf
Tangent vectors to normal slices are the values of those Jacobi
fields along normal geodesics vanishing on $N$ and with
derivative there in $\nu N.$  Let $X(t)$ be such a field
along the geodesic $\gamma(t)$, parametrized by arc length.  The
initial derivative $X'(0)$ belongs to $\nu_{\gamma(0)}N$; for
convenience we take it to be a unit vector.  As in the proof of
Proposition \ref{prop-focal}, we identify all the tangent spaces
along $\gamma$ with $T_{\gamma(0)}M$ by parallel translation.

Since $X(0)=0$, we write $X(t)=tX'(0)+Z(t)$.  Then $Z(0)=Z'(0)=0$,
and $Z''(t)=X''(t)=-R(t)X(t)$ by the Jacobi equation.  The
curvature estimate $\|R(t)\| \leq 1$ gives first of all the
estimate $\|X(t)\| \leq \sinh t$, so $\|Z''(t)\|\leq \sinh t$.
Integration gives $\|Z'(t)\|\leq \cosh t -1$ and then
$\|Z(t)\|\leq \sinh t - t \leq t^3/6 $.  Then $\|X(t)/t -
X'(0)\|\leq t^2/6$ and $\|X(t)/t\|\geq 1-t^2/6$.

The typical quasi-vertical unit vector $(X(t)/t)/\|X(t)/t\|$ is thus
within distance $$\frac{t^2/6}{1-t^2/6} \leq \frac {t^2}{5}$$ of the
vertical vector $X'(0)/\|X(t)/t\|$, from which it follows (for $t\leq
1$) that the distance between the vertical and quasi-vertical spaces
at $x=\gamma(t)$ is at most
$$ \sin^{-1}(t^2/5) \leq t^2/4=d(x,N)^2/4,$$ while the
difference between the projections on these spaces is at most
$\frac{1}{5}d(x,N)^2$.
\qed

An estimate for the covariant derivative of the quasi-vertical bundle
in non-vertical directions would simplify our later considerations, but
we have not been able to find such an estimate.\footnote{The
difficulty seems to be related to the fact that a bound on curvature
does not imply a bound on the derivatives of Christoffel symbols in
normal coordinates; see the example on p.34 of
\cite{jo-ka:geometrische}.}
Instead, we will
estimate the covariant derivatives of the vertical bundle and then use
Proposition \ref{prop-twobundles} to shuttle between the two bundles.

\begin{prop}
\label{prop-covderproj}
Let $(M,N)$ be a gentle pair, $\Gamma_N$ its extended Gauss map, and
$\proj_{\Gamma_N}$ the field of orthogonal projections onto the image
of $\Gamma_N$ (the vertical bundle).  Then the covariant derivative
$\covder\proj_{\Gamma_N},$ a field of maps from $TM$ to ${\rm
End}~TM$, has operator norm bounded pointwise by $11$ on the tube
$\calt^{\frac{1}{2}}N.$
\end{prop}
\pf
We estimate $\covder \proj_{\Gamma_N}$ along $N$ in terms of the second
fundamental form of $N$, then we extend the estimate into the tubular
neighborhood by differentiating along the normal geodesics.

Along $N$, $\covder_v\proj_{\Gamma_N}$ is zero if $v$ is vertical, by
the definition of $\Gamma_N$.  If $w$ is horizontal,
$\covder_w\proj_{\Gamma_N}$ is a symmetric operator which exchanges
horizontal and vertical spaces.  For vertical $v$ and horizontal $h$,
it follows easily from the definition of the second fundamental form
that $\langle (\covder_w\proj_{\Gamma_N}) (v),h\rangle
= B_v(w,h).$  Along $N$, then, $\|\covder\proj_{\Gamma_N}\|=\|B\|$,
which is at most $\frac{3}{2}$ by Corollary \ref{cor-secondbound}.

To estimate $\covder\proj_{\Gamma_N}$ at a point $x\in \calt N,$ we
let $\gamma$ be the unit speed geodesic segment from $\pi_N(x)$ to $x$
and estimate $\covder_{\gamma ' (s)} \covder \proj_{\Gamma_N}.$  Let
$Y$ be a unit length vector field on a neighborhood of $\gamma$ such that
$\gamma$
is an integral curve of $Y$.   For any unit vector $v\in T_x M$, we can
find a vector field $Z$ near $\gamma$ such that $[Y,Z]=0$ and
$Z(x)=v$.

Now $\covder_Y \covder_Z \proj_{\Gamma_N} = \covder_Z
\covder_Y \proj_{\Gamma_N} + [R(Y,Z), \proj_{\Gamma_N}].$  But
$\covder_Y  \proj_{\Gamma_N} = 0,$\linebreak $\|Y\|=1$, and $\|
\proj_{\Gamma_N}\|=1$, so $\|\covder_Y \covder_Z  \proj_{\Gamma_N} \|
\leq 2 \|R\|~\|Z\|$, and hence
%\begin{equation}
%\label{eq-covest1}
$$\|\covder_v  \proj_{\Gamma_N}\| \leq \textstyle{\frac{3}{2}} +  d(x,N)
\cdot 2\|R\| \displaystyle{\sup_{0\leq
  s \leq d(x,N)}} \|Z(\gamma(s))\|.
$$
%\end{equation}

Since the sectional curvatures are bounded in absolute value by 1, \
$\|R\|\!\leq\!\frac{4}{3}$ (see \S6.1 in \cite{bu-ka:gromov}), so it
remains only to estimate the size of $Z$ along $\gamma$.  For this
purpose, we construct a particular choice of $Z$ by forming the
rectangle (as used in Lemma \ref{lem-warner}) consisting of minimizing
geodesics from $N$ to the points on the geodesic with initial vector
$Z(x)=v$ and letting $Z$ be the variational vector field.  $Z$ is an
$N$-Jacobi field along $\gamma$, so we can estimate it using Warner's
comparison theorem \cite{wa:extensions}.
As in Lemma \ref{lem-warner}, we conclude from Warner's
theorem that $\|Z(\gamma(s))\| \geq (\cos s - \frac{3}{2} \sin s)
\|Z(\gamma (0))\|$ for all $s$.  Setting $s=d(x,N)$, we get:
$ \|Z(\gamma(0))\| \leq (\cos d(x,N) -\frac{3}{2} \sin d(x,N))^{-1},$
which is less than $7$ if $d(x,N) < \frac{1}{2}$.

Now the maximum of $\|Z(\gamma(s))\|$ might occur at an interior point
$\gamma(s_0)$ of the segment.  If that happens, we must have
$\covder_{\gamma'(s_0)} Z \perp Z(s_0)$.  In that case, we can use
Warner's theorem again, comparing with the situation where the initial
manifold is totally geodesic, to get an estimate even better than
before.

Our conclusion is that $ \sup_{0\leq
  s \leq d(x,N)} \|Z(\gamma(s))\| \leq 7 $, so that
$\|\covder  \proj_{\Gamma_N} \| \leq \frac{3}{2} + \frac{1}{2}\cdot 2
\cdot \frac{4}{3} \cdot 7 < 11.$
\qed

\subsection{The potential function of a submanifold}

On the tubular neighborhood
$\calt N$, the distance function from $N$,
$\rho_N(x)=d(x,\pi_N(x))=\|\exp_N^{-1}(x)\|,$ is singular along $N$,
but $P_N = \frac{1}{2} \rho_N^2$ is smooth
and has $N$ as a nondegenerate critical manifold (of absolute minima).
We call $P_N$ the {\bf potential function} of $N$, thinking of $N$ as
the attracting manifold of a ``linear'' attracting force field.
The points of $N$ can also be recovered from $P_N$ as the
nondegenerate minimum points of its restrictions to the normal slices,
to each normal slice, or equivalently as the zeros of the projection
of the gradient $\grad P_N$ into the quasi-vertical bundle.  This
description of $N$ will be most useful to us, since the projected gradient is
transverse to the zero section in the quasi-vertical bundle.  We will be
proving Theorem \ref{thm-centerofmass} by identifying the invariant
submanifold $\overline{N}$ as the zero set of an invariant object
obtained by averaging $\grad P_N$ over $G$ and projecting into an
averaged vertical bundle.

We will need estimates on the hessian $H_N=\covder \grad P_N$, a field
of symmetric bilinear forms on the tangent spaces of $\calt N$ defined by
$H_N(X,Y)=\langle \covder_X\grad P_N,Y\rangle.$
 Estimates for this hessian when $N$ is a
 point are given in \cite{bu-ka:gromov} (Proposition 6.4.6)
   and play an important role in applications
of the center of mass construction. When $\dim N \geq 1$, the
estimates are more complicated, since the hessian behaves quite
differently in horizontal and vertical directions.

To calculate the hessian of $P_N$, we adapt the method of Gray
\cite{gr:tubes}, who used a Riccati equation to study the second
fundamental forms of the boundaries of the $r$-tubes, i.e. of the
level hypersurfaces of $P_N$.  His analysis used the gradient $U_N$ of
the function $\rho_N$, which has the same level hypersurfaces.
Since $U_N$ is the unit outward normal to these hypersurfaces, the
hessian $S_N=\covder U_N$ gives the second fundamental form of the tube
boundaries.  The covariant derivative of the tensor field $S_N$ in the
radial direction satisfies the Riccati equation
\[\covder_{U_N}S_N=S_N^2+R_N,\]
 where $R_N$ is the field of operators obtained (as one does
for the Jacobi equation)
by contracting the radial
vector field twice with the Riemann
curvature tensor $R$ of $M$; i.e.
$R_N(X)=R(U_N,X)U_N.$ The Riccati equation determines $S_N$ in the
tubular neighborhood once we know its behavior along $N$.  But $S_N$,
like $U_N$, is not even defined along $N$, so this initial behavior
must be described asymptotically.  This point is dealt with only
briefly by Gray, so we paraphrase here the following result of
Eschenburg (\cite{es:comparison}, \S6.1).

\begin{prop}
\label{prop-eschenburg}
Along a unit speed geodesic $\gamma$ emanating normally from $N$, the
solution $S_N$ of the Riccati equation has the asymptotic expansion
$S_N(\gamma(t))=t^{-1} S_{-1} + S_0 + {\rm O}(t)$, where $S_{-1}$ and
$S_0$ are in block matrix form with respect to the decomposition into
horizontal and vertical subspaces as follows:
$$
S_{-1} =
\left[\begin{array}{c|c}  I & 0\\ \hline 0 & 0
\end{array}\right]
~~~~~~~~
S_0=
\left[\begin{array}{c|c} 0  & 0\\ \hline 0 & B_{\gamma'(0)}
\end{array}\right].
$$

\end{prop}

Unfortunately, the proof of this Proposition
seems to belong to folklore rather than literature.  Nevertheless,
we will not write a proof here, either.  (The proof is a nice
exercise in the use of Fermi coordinates or Fermi frames \cite{gr:tubes}.)

We return now to
the hessian $H_N=\covder \grad P_N$.  Differentiating the relation
$\grad P_N=\rho_N \grad \rho_N$ yields
$$H_N=U_N\otimes U_N + \rho_N S_N,$$
 where $U_N\otimes U_N $ is the
operator of orthogonal projection onto the unit vector field $U_N$.
As the hessian of a smooth function, $H_N$ is defined throughout the
tubular neighborhood, and its values along $N$ are quite easy to
compute.  In fact, since $P_N$ is critical at all points of $N$, we
can compute the hessian at any such point as the second derivative
matrix in any coordinate system for which the coordinate vector fields
form an orthonormal basis at that point.  From this observation, it is
easy to check (see below) that $H_N$ is block diagonal with respect to
the normal--tangent splitting of $TM$ along $N$.  On the tangent
bundle $TN$, $H_N$ is zero, while it equals the riemannian metric
on the normal bundle $\nu N$.

Since it is not $H_N$ but rather the singular field
$$S_N = (1/\rho_N) ( H_N-U_N\otimes U_N)$$
which satisfies the Riccati
equation, the initial conditions required for determining or
estimating a solution also include the normal (covariant) derivative of $H_N$
along $N$.  It is here that the second
fundamental form of $N$ appears, as we saw in Proposition
\ref{prop-eschenburg} and  one may easily check in the
simple example where $N$ is a circle in the euclidean plane.

\subsection{Estimates for the hessian}
We will use the following estimate of Eschenburg and Heintze
\cite{es-he:comparison} on solutions of Riccati equations to estimate
the hessian of the potential function.
\begin{thm}[Eschenburg--Heintze \cite{es-he:comparison}]
\label{thm-riccati}
Let $S_1(t)$ and $S_2(t)$ be\newline maximal solutions of the Riccati
equations
$$S_j'=S_j^2+R_j$$
 defined for $0<t<t_j$ with values in symmetric
matrices such that:\newline
(1) U := $S_2-S_1$ has a continuous extension to
$t=0$ with $U(0)\geq 0$;\newline (2) $R_j(t)$ is a smooth function with
$R_2(t)\geq R_1(t)$ for all $t\geq 0$.\newline  Then $t_1 \leq t_2$, and $S_1
\leq S_2$ on $(0,t_1)$.
\end{thm}
We recall that an inequality of the form $P\leq Q$ between symmetric
matrices means that the difference $Q-P$ is positive semidefinite.  It
implies that the maximal eigenvalue of $P$ is less than that of
$Q$. As a consequence, the inequality $-cI\leq P \leq cI$ for a
positive number $c$ implies the bound $\|P\|\leq c$ on the operator
norm of $P$.

We will use Theorem \ref{thm-riccati} to estimate the hessian of the
potential function for a gentle pair,
by comparison with manifolds of constant curvature $\pm 1$.  Our
first step will be to examine and estimate the solutions of the
Riccati equation when $R(t)$ is a constant multiple of the identity
matrix.

It suffices to consider the scalar case, i.e. the equations
$$s'=s^2\pm 1$$
 for a real-valued function $s(t)$ defined for $t>0$.
We must consider two kinds of boundary conditions, corresponding to
the vertical and horizontal cases in the geometric situation.

For the vertical case (directions normal to $N$ and their parallel
translates along normal geodesics), we choose the solution with a
pole at 0, namely $-\cot t$ when the sign of the curvature $\pm 1$
is positive and $-  \coth  t$ when it is negative.  For convenience,
we will sometimes use the notation $\cot_+$ for $\cot$ and $\cot_-$
for $\coth$, as well as $\cot_\pm$ to denote whichever is appropriate;
thus the general solution can be written as
\[s_{\ver}(t)= \cot_{\pm} t.\]

For the horizontal case (directions tangent to $N$ and their
parallel translates along normal geodesics),
we take the
 solution with initial value $b$
(which can be positive or negative), which is
\[s_{\hor}(t)=\frac { \tan_{\pm} t+b}{1\mp b\tan_{\pm} t}.\]

We will now estimate the function $h(t)=ts(t)$ on the interval
$[0,1 ]$, since it is this product rather than $s(t)$ alone which
gives the hessian of the potential function.

For the vertical direction, we have $h_{\ver}(t)= t\cot_{\pm} t$
(notice that the pole has disappeared, as it should).  With positive
curvature, $h(t)$ decreases from $1$ at $ t=0$ to $\cot 1 \sim
.642\ldots$ at $ t=1.$ With negative curvature, $h(t)$ increases from
$1$ at $ t=0$ to $\coth 1\sim 1.313\ldots$ at $ t=1$.  Thus, we have
the estimate
\[.64<h_{\ver}(t)<1.32~~{\rm f}{\rm o}{\rm r}~~0\leq t\leq 1 .\]
 Notice
that this estimate also applies to the hessian in the radial
direction, which is identically $1$, even though the relation between
the hessian and the Riccati solution is more complicated in that case.

For the horizontal direction, we will need only an upper bound, and
there is no pole to worry about, so we can estimate $s_{\hor}$
first
and multiply by $t$ afterward.  For positive curvature,
\[s_{\hor}(t)=\frac { \tan  t+b}{1-b\tan  t}\]
 is
monotonically increasing in $b$, so we can set $b=\frac{3}{2}$, by our
estimate on the second fundamental form under the assumption of
bounded geometry.  This gives us
\[|s_{\hor}(t)|< \frac {\tan  t+\frac{3}{2}}{1-\frac{3}{2}\tan t} \ .\]
 To
keep the denominator under control, we will require that $\frac{3}{2}
\tan t
< \frac{1}{2}$, i.e. $\tan t < \frac{1}{3}$, which is satisfied when $ t <
.321\ldots$.  To keep the numbers simple, we will require that $ t <
\frac{1}{4}$, which leads to the upper bound of $\frac{13}{6}$, which for
simplicity
we replace by $3 $.  This gives us the estimate
\[|h_{\hor}(t)|\leq 3 t~~{\rm f}{\rm o}{\rm r}~~0\leq t\leq
\textstyle{\frac{1}{4}} \ .\]

We are left with the negative case.  But since $\tanh$ grows more
slowly than $\tan$, the estimate from the positive case still applies.

Summarizing the estimates above and applying Theorem
\ref{thm-riccati}, with the initial conditions given by Proposition
\ref{prop-eschenburg}
gives the following:

\begin{prop}
\label{prop-hestimate}
Let $(M,N)$ be a gentle pair.  At a
point of $M$ with distance $r\leq \frac{1}{4}$ from $N$, the hessian $H_N$
when expressed in block matrix form with respect to the decomposition
into vertical and horizontal parts satisfies the
inequality:
\begin{equation}
\left[\begin{array}{c|c} .64 I & 0\\ \hline 0 & -3 rI
\end{array}\right]
< H_N < \left[\begin{array}{c|c} 1.32 I & 0\\ \hline 0 & 3 rI
\end{array}\right]
\end{equation}
\end{prop}
\begin{cor}
\label{cor-hessian}
On the $\frac{1}{4}$-tube around $N$, $\|H_N\|\leq 1.32$.

\end{cor}
\pf The statement follows from the bound $-.75I \leq H_N \leq 1.32 I$.
\qed

\begin{cor}
\label{cor-hv}
At a point of $M$ with distance $r\leq\frac{1}{4}$ from $N$, for a
horizontal vector $v$ and a vertical vector $w$, $|H_N(v,w)| \leq 3 \sqrt{r}
~\|v\|~\|w\|.$
\end{cor}
\pf
We may assume that $v$ and $w$ are unit vectors.  By Proposition
\ref{prop-hestimate}, $|H_N(v,v)| < 3r$ and $ .64 < H_N(w,w) < 1.32$.  In
addition, for any real number $t$, $H_N(tv+w,tv+w) < 3rt^2 + 1.32$.
Thus the quadratic function $t\mapsto
t^2 (H_N(v,v)-3r)  +2tH_N(v,w) + H_N(w,w) - 1.32 $ is negative for all
$t$ and hence its discriminant is negative.  So
$4H_N(v,w)^2 < 4(H_N (v,v) -3r) (H_N (w,w) -1.32) \leq 4\cdot 6r\cdot 0.68.$
Therefore, $H_N^2(v,w) < 4.08 r$, and so $|H_N(v,w)| < 3 \sqrt r$.
\qed

\begin{rmk}
\label{rmk-squareroot}
{\em
The square root in Corollary \ref{cor-hv} is responsible
for the presence of $\sqrt{\epsilon}$ rather than $\epsilon$ in the
estimates in the main theorems.   It seems plausible that the square
root could be removed here, since the off-diagonal block of $H_N$ is
zero along $N$, and so it must grow linearly with $r$ on any particular
manifold.  Because of the singular nature of the relation between
$S_N$ and $H_N$, though, we have not been able to get a linear bound which
is uniform for all gentle pairs.
}
\end{rmk}

\subsection{Geometry relative to nearby submanifolds}

In this section, we estimate the difference between geometric
objects in the tubular neighborhoods of two nearby submanifolds.
To begin, we show that the retractions to nearby submanifolds are
close to one another.

\begin{lemma}
\label{lem-polyestimate}
Let $(M,N)$ and $(M,N')$ be gentle pairs such that $d(N,N')<\frac{1}{2}$,
and let $x\in M$ with $d(x,N) <\frac{1}{4} $.
 Then $d(\pi_N(x),\pi_N \pi_N'(x)) <14 d(N,N')$ and
  $d(\pi_N(x),\pi_N'(x)) <15 d(N,N')$.
\end{lemma}
\pf Let $w=\pi_N(x)$, $y=\pi_{N'}(x)$, and $z=\pi_N(y)$.
Let $\xi$ be the unit tangent
vector at $y$ to the geodesic segment $yx$, let $\eta$ be the
parallel translate of $\xi$ to $z$ along the geodesic segment
$yw$, and let $\zeta$ be the nearest unit vector to $\eta$ in the normal space
$\nu_z N$.  Let $r$ and $s$ be the endpoints of the geodesic segments
$zr$ and $zs$ with initial directions $\eta$ and $\zeta$ respectively,
each with the same length as $yx$.  We will show below that
$d(x,r)$ and $d(r,s)$ are small.  The resulting estimate on
$d(x,s)$, combined with Lemma \ref{lem-warner}, will give us an
estimate on $d(w,z)$, which will lead to the result.

To estimate $d(x,r)$, we observe that $x$ and $r$ are the
endpoints of geodesic segments departing from $y$ and $z$ with
initial tangent vectors related by parallel translation along the
geodesic $yz$.  It follows that $d(x,r)\leq \lambda d(y,z)$, where
$\lambda$ is a bound on the growth of a Jacobi field with initial
derivative 0 along a
geodesic segment of length $d(x,y)$,
under the assumption that sectional curvatures are $\geq -1$.

We write $\alpha$ for $d(N,N')$ and $\beta$ for
$d(x,N)$.  Since $d(x,y)\leq d(x,N) + d(N,N') = \beta + \alpha$,
we can take $\lambda = \cosh (\beta+\alpha)$, and since
$d(y,z)\leq d(N,N')=\alpha$ we have
$$d(x,r) < \alpha \cosh(\beta+\alpha).$$

For $d(r,s)$, we use the fact that $r$ and $s$ are the endpoints
of geodesic segments with the same initial point, length at most
$\alpha + \beta$, and initial (unit) tangent vectors making an angle of at
most $\alpha$ with one another (since $d(\Gamma_N (z),\Gamma_{N'} (z))\leq
d(N,N')=\alpha$).   The usual Jacobi field estimate, using again
the lower bound of $-1$ on curvature, gives $$d(r,s) < \alpha \sinh
(\beta + \alpha),$$
and hence
$$d(x,s) < \alpha(\cosh(\beta+\alpha) + \sinh(\beta +\alpha))
=\alpha e^{\beta+\alpha} <\alpha e^{\frac{1}{3}} < 2 \alpha.$$

By Corollary \ref{cor-seven}, $d(w,z) < 14\alpha$.  Since
$d(y,z)\leq\alpha$, the triangle inequality gives the required
estimate for $d(\pi_N(x),\pi_N'(x))$.
\qed
Next we estimate the difference between the tangent vectors from a
point $x$ to the submanifolds $N$ and $N'$.
\begin{lemma}
\label{lem-gradientdifference}
Let $(M,N)$ and $(M,N')$ be gentle pairs such that $d(N,N')<\frac{1}{50}$,
and let $x\in \calt N$ be such that $d(x,N)<\frac{1}{4} .$
Then the difference
between the initial tangent vectors
 $\xi_N(x)$ and $\xi_{N'}(x)$  of the minimizing geodesics from
 $x$ to $\pi_N(x)$ and $\pi_{N'}(x)$ respectively is less than
$\frac{45}{2} ~d(N,N')$.
\end{lemma}
\pf
As in Lemma \ref{lem-polyestimate}, we let $\alpha=d(N,N')$.
According to that lemma, we can connect $\pi_N(x)$ to
$\pi_{N'}(x)$ by a path $\sigma$ of length less than $15\alpha <
\frac{3}{10}$.  Since the
distance from $x$ to $N$ is less than $\frac{1}{4}$, all the points
of $\sigma$ are within a distance $1$ of $x$, so that we can lift
$\sigma$ by the exponential map at $x$ to a path $\tilde{\sigma}$
in the tangent space $T_x M$ and thereby produce a smooth family
$\{\gamma_s\}$ of geodesic segments joining $x$ to the points of
$\sigma$.  Note that  $\tilde{\sigma}(0)$ and $\tilde{\sigma}(1)$
equal  $\xi_N(x)$ and $\xi_{N'}(x)$ respectively.

We estimate the length ratio $|\sigma|/|\tilde{\sigma}|$ by using the
Jacobi fields along the segments $\gamma_s$.  Since the sectional
curvature in $M$ is bounded above by 1, the segments have
length less than 1, and $\sin x/x > \frac{2}{3}$ for $0\leq x\leq 1$, we
have the estimate
$$|\sigma|/|\tilde{\sigma}|
> \textstyle{\frac{2}{3}}.$$   $\frac{3}{2} \cdot 15 \alpha =
\frac{45}{2}\alpha$ is therefore an upper bound for the
length of $\tilde{\sigma}$ and hence for the difference between
$\xi_N(x)$ and $\xi_{N'}(x).$
\qed

Finally, we estimate the distance between the extended Gauss maps
of nearby submanifolds.

\begin{prop}
\label{prop-verticalangle}
Let $(M,N)$ and $(M,N')$ be gentle pairs such that \linebreak
$d(N,N')\!\leq~\frac{1}{2}$. If $x \in \calt_{\frac{1}{4}}N$,
then $d(\Gamma_N(x),\Gamma_{N'}(x)) < 250 d(N,N').$
\end{prop}
\pf
Let $y=\pi_{N'}(x)$, $w=\pi_N (x)$, and $z=\pi_N(y)=\pi_N \pi_{N'}(x)$.
 The polygon $xyzwx$ is defined to consist of segments which
are geodesics in $M$, except for $zw$, which is taken to be a
geodesic in $N$.
 We lift the
polygon $xyzwx$ to $GM$ as follows.  On the segments $xy$ and $yz$, we take the
parallel section $\Gamma_{N'}$; note that
$$d(\Gamma_{N'}(z),\Gamma_N(z))\leq d(N,N')<\alpha.$$
 At this point,
we jump to $\Gamma_{N}(z)$, which we parallel
translate to $w$ along the geodesic $zw$ in $N$.  Arriving at $w$,
we have a space whose distance from $\Gamma_{N}(w)$ is
less than $\frac{3}{2}$ (our bound for the second fundamental form of
$N'$) times the length $|zw|$.  Finally, we jump to $\Gamma_{N}(w)$
and parallel translate it to $x$ along the normal geodesic $wx$ to
arrive at $\Gamma_{N}(x)$.  The distance between the
starting and ending values of our section is bounded above by
the sum of the sizes of the jumps at $z$ and $w$ and $2^{-\frac{1}{2}}$ times
the distance from the identity to the operator $\calh$
of holonomy around $xyzwx$.

Let $\alpha=d(N,N')$. By  Lemma \ref{lem-polyestimate},
$d(z,w) < 14\alpha$; hence, $|zw|
<7\cdot14\alpha <100\alpha$ by Corollary \ref{cor-mndistance}, and so
the jump at $w$ has size at most $150\alpha$.  The distance
$d(I,\calh)$ from the identity to $\calh$ is
bounded by the product of the curvature ($\leq 1$) and the area of a
surface spanning the polygon.  We build such a surface\footnote{This is one
of the places where we use Condition (iii) in the definition of bounded
geometry; see Remark \ref{rmk-condition}} with the
minimizing geodesics from $x$ to the points on the segments $yz$ and
$zw$.  All these points are within a distance $d(x,N)+2\alpha
<\frac{1}{4}+2\alpha < \frac{3}{4}$ of $x$, and
Jacobi fields along these geodesic segments with the value zero at $x$
are always increasing in length (no focal points), so the area of the
 surface is at most
$$\textstyle{\frac{1}{2}}(d(x,N)+2\alpha)(|yz|+|zw|)<
(\textstyle{\frac{1}{4}}+1)(\alpha+100\alpha)<
130\alpha. $$

Adding the estimates for the jumps and the holonomy, we conclude that
$$d(\Gamma_{N}(x),\Gamma_{N'}(x))<\alpha+150\alpha +2^{-\frac{1}{2}}\cdot
130\alpha
<250\alpha .$$
\qed

\begin{rmk}
\label{rmk-quasisymmetry}
{\em
We can derive from the results of this section a crude\linebreak
``quasisymmetry'' estimate for the $C^1$ distance: if  \ $d(N,N') <
\frac{1}{4}$, then\linebreak
$d(N',N) < 250d(N,N')$.
}
\end{rmk}

\section{The main theorem}
\label{sec-main}
In this section, we will prove Theorem \ref{thm-centerofmass}, of
which Theorem \ref{thm-main} is a direct consequence.  For convenience, we
repeat
the statement of the former:

\bigskip
\noindent
{\bf Theorem \ref{thm-centerofmass}}~
{\em Let $M$ be a riemannian manifold and $\{N_g\}$ a family of compact
submanifolds of $M$ parametrized in a measurable way
by elements $g$ of a probability space $G$, such that all
the pairs $(M,N_g)$ are gentle.  If $d(N_g,N_h) <
\epsilon < \frac{1}{20000}$ for all $g$ and $h$ in $G$, there is a well
defined {\bf
center of mass} submanifold $\overline N$ with $d(N_g,\overline {N}) <
136\sqrt{\epsilon}$ for all $g$ in $G$.
The center of mass construction is equivariant with respect to
isometries of $M$ and measure preserving automorphisms of $G$.
}

\bigskip
\pf
We begin with an outline of our proof and leave details
to a series of lemmas which follow.

For each $g$ in the parameter space $G$, the potential function
$P_g=P_{N_g}$ is defined on the tubular neighborhood $ \calt
N_g$.  On the intersection $\calt$ of these
neighborhoods, the average over $g$ of the functions $P_g$
is a function $\overline{P}$ whose gradient is the average of
the $\grad P_g$.  The Gauss maps $\Gamma_g=\Gamma_{N_g}$
are all defined on $\calt$ as well, and since the manifolds $N_g$ are $C^1$
close to one another, the average of the $\Gamma_g$ (by the
construction in the Appendix) is a
well-defined section $\overline{\Gamma}$ of the Grassmann bundle.
The zero set $\overline{N}$ of the projection
$\proj_{\overline{\Gamma}}\grad \overline {P}$ of
$\grad \overline{P}$ into the averaged vertical bundle $\overline{\Gamma}
(\calt)$ is
obviously an invariant subset of $M$.

To see that $\overline{N}$ is a smooth manifold which is
$C^1$ close to all the $N_g$, we choose a base point $e$ in $G$ and
look at the restriction of $\proj_{\overline{\Gamma}}\grad
\overline{P}$ to each of the normal
slices for $N_e$.  To study these restrictions, we make
$\proj_{\overline{\Gamma}}\grad
\overline{P}$
 into a vector
field along the normal slices by applying the projection $\proj_{Q_e}$ into
the quasi-vertical bundle (tangents to the normal slices) for $N_e$.
Since the vertical and quasi-vertical bundles are close, by
Proposition
\ref{prop-twobundles}, $\proj_{Q_e}$ gives an isomorphism between
the two bundles, and $\overline{N}$ is
again the zero set of the vector field
$\xi=\proj_{Q_e}\proj_{\overline{\Gamma}}\grad
\overline{P}$
along the normal slices.  (We note that $\xi$ is
neither invariant nor a gradient, but this is not a problem.)

Combining several estimates from Section \ref{sec-geometry}, we will
show
that all the zeros of $\proj_{Q_e} \proj_{\overline{\Gamma}}
\overline{\grad P}$
are nonsingular with
index $+1$ along the normal slices, and that
$\proj_{Q_e} \proj_{\overline{\Gamma}}\grad \overline{P}$
points outward along the sphere of radius
$100 \epsilon$ in each slice.  As a result, there is a unique zero inside that
sphere, and the collection of all these zeros forms a manifold which
is transverse to the normal slices and, in fact, is $C^1$ close to
$N_e$.   Since $e$ was arbitrary, we get $C^1$ closeness of $\overline{N}$ to
all the $N_g$, and we are done.

\begin{rmk}
\label{rmk-infinite}
{\em
The arguments in the preceding paragraph require $M$ to be finite
dimensional.  To extend the result to the infinite-dimensional case
will require further analytic estimates on the projected gradient, so
that some infinite-dimensional degree theory can be applied.
}
\end{rmk}

We will divide the detailed proof into a series of steps, using
the assumption that the $(M,N_g)$ are gentle pairs with
$d(N_g,N_h) < \epsilon$ for all $g$ and $h$.  We will add restrictions
on $\epsilon$ as necessary.  Lemmas will be embedded in the proof;
the end of the proof of a lemma will be indicated by the symbol $\triangle$.

\bigskip
\noindent
{\bf STEP 1.} To begin, we establish a region on which all the
potential functions are defined.
\begin{lemma}
\label{lem-regions}
For any $r$ such that $\epsilon < r < 1$, $\calt ^r =
\cap_{g \in G} \calt^r N_g$ contains $\calt ^{r-\epsilon} N_h$ for
every $h$ in $G$.
\end{lemma}
\pf
If $x\in \calt ^{r-\epsilon} N_h$, then there
 is a $y \in N_h$ with $d(x,y)
< r-\epsilon$.  For any $g$, since $d(N_g,N_h)<\epsilon$, there is a $z$ in
$N_g$
with $d(z,y) < \epsilon$.  Then $d(x,z) < r$, and so $x \in \calt^r
N_g$.  Thus $\calt ^{r-\epsilon} N_h$ is contained in each $ \calt^r N_g$.
\lqed
We will leave $r$ undetermined for the moment and study geometry in
the $G$-invariant subset $\calt^r$.

\bigskip
\noindent
{\bf STEP 2.}  The Gauss maps
$\Gamma_g=\Gamma N_g$ are all defined on $\calt^r$.  We will show
that, when $r$ and $\epsilon$ are
 sufficiently small, they are close enough to one another
so that their values can be averaged by the method of the Appendix to
produce an equivariant Gauss map $\gammabar$.  This averaged Gauss map
is close to both the vertical and quasi-vertical bundles of the
$N_g$.
\begin{lemma}
\label{lem-closeforgauss}
If $r < \frac{1}{4}$ and $\epsilon < \frac{1}{1000}$, then the average
$\gammabar(x)$ of the
subspaces $\Gamma_g(x) \subset T_x M$ can be constructed by
Proposition \ref{prop-averageproj} for each $x \in \calt^r$, with
the Finsler distance $d(\gammabar(x), \Gamma_g(x)) < 1000 \epsilon$ for all
$g$ in
$G$.  The mapping $\gammabar : \calt^r \rightarrow TM$ is smooth and
$G$-equivariant, and the construction of $\gammabar$ from the family
$\{N_g\}$ is equivariant with respect to isometries of $M$ and measure
preserving automorphisms of $G$.
\end{lemma}
\pf
Since $x$ has distance less than $\frac{1}{4}$ from
 each $N_g$, and $d(N_g,N_h) < \epsilon$,
Proposition \ref{prop-verticalangle}  gives the estimate
$d(\Gamma_g(x),\Gamma_h(x)) < 250 \epsilon < \frac{1}{4}$.
 Proposition \ref{prop-averageproj} then produces in each tangent space
$T_x M$ an averaged subspace $\gammabar(x)$ with $d(\gammabar(x) ,
\Gamma_g(x)) < \sin^{-1} (3\cdot 250\epsilon)$ for all $g$.  Since
$\sin^{-1} x/x < \frac{4}{3}$ when $0 < x < \frac{3}{4}$, we conclude that
 $d(\gammabar(x) ,\Gamma_g(x)) < 1000 \epsilon$.

The equivariance properties of $\gammabar$ and its construction
follow from the absence of arbitrary choices in the construction.  The
smoothness of $\gammabar$ follows from the corresponding assertion in
Proposition \ref{prop-averageproj} (using local trivializations).
\lqed
\begin{lemma}
\label{lem-quasiandaverage}
For each $x$ in $\calt^r$ and $g$ in $G$, the distance from
$\gammabar(x)$ to the quasi-vertical space $Q_g=\ker T_x \pi_{N_g}$ is
less than $1000 \epsilon + \frac{1}{4}r^2$.
\end{lemma}
\pf
Add the estimates in Propositions \ref{prop-twobundles} and
\ref{prop-averageproj}.
\lqed

\bigskip
\noindent
{\bf STEP 3.}  The averaged potential function
$$\pbar = \int_G P_g~dg$$ is defined on $\calt^r$.  Since the metric
is invariant, $$\grad \pbar = \int_G \grad P_g ~dg.$$
\begin{lemma}
\label{lem-potentials}
If $r < \frac{7}{10}$, $\epsilon < \frac{1}{50}$, and $x\in \calt^r$,
then $$\|\grad \pbar (x)- \grad P_g (x)\| < \frac{45}{2} \epsilon. $$
\end{lemma}
\pf
By the usual averaging argument, it suffices to check that $$
\|\grad P_g - \grad P_h\| < \frac{45}{2} \epsilon$$ for all $g$ and
$h$ in $g$.  Observing that $-\grad P_g(x)$ is just the initial tangent
vector of the minimizing geodesic from $x$ to $N_g$, we get the result
immediately from Lemma \ref{lem-gradientdifference}.
\lqed

\bigskip
\noindent
{\bf STEP 4.}
As indicated in our outline of the proof, we define $\nbar$ to be the
set of zeros in $\calt^r$ of the projected gradient vector field
$\proj_{\overline{\Gamma}}\grad \overline {P}.$   It is obviously a
$G$-invariant subset of $M$.  To study this zero set, we fix a
reference element $e$ of $G$ and replace
$\proj_{\overline{\Gamma}}\grad \overline {P}$ by its projection
$\proj_{Q_e}\proj_{\overline{\Gamma}}\grad \overline {P}$ into the
quasi-vertical bundle $Q_e$ for $N_e$.  Since this projected field
is tangent to the normal slices for $N_e$, it will have at least
one zero in each normal slice if we can show that it points
outward along some sphere, i.e. if it has a positive inner product
with $\grad P_e$.
As long as the hypotheses of Lemma \ref{lem-closeforgauss} are
satisfied, Lemma \ref{lem-quasiandaverage} implies that the distance
$\|\proj_{\overline{\Gamma}}-\proj_{Q_e}\|$ is less than $\sin (1+\frac{1}{16})
< 1$, so the projection from the averaged normal bundle into the
normal slices is an isomorphism, and the zero set is not changed.  To
simply formulas, we will denote the ``doubly projected gradient''
vector field $ \proj_{Q_e}\proj_{\overline{\Gamma}} \grad \pbar
(x)$ by the symbol $\doub$.

\begin{lemma}
\label{lem-innerproduct}
If $\epsilon<\frac{1}{2250}$ and $99\epsilon< d(x,N_e) < \frac{1}{4},$
then the inner product
$\langle \doub,
    \grad P_e(x) \rangle $ is positive.
\end{lemma}
\pf
Let $\rho=d(x,N_e)$.  Under the assumptions of the Lemma, it
follows from Lemma \ref{lem-quasiandaverage} and Proposition
\ref{prop-finslervsnorm} that
$\|\proj_{Q_e} -\proj_{\gammabar}\| < \sin (\frac{1000}{2250} +
\frac{1}{16}) < \frac{1}{2},$ and from Lemma
\ref{lem-gradientdifference} that $\|\grad \pbar - \grad P_e\| <
\frac{45}{2}\epsilon < \frac{45}{198}\rho.$

Then
\begin{eqnarray*}
\langle \doub,
    \grad P_e(x) \rangle &=& \langle
\proj_{Q_e}\proj_{\overline{\Gamma}} \grad \pbar (x),
    \grad P_e(x) \rangle   =
 \langle  \grad \pbar (x),
  \proj_{\overline{\Gamma}} \grad P_e(x) \rangle \\ & = &
  \langle  \grad P_e(x) + (\grad \pbar (x) - \grad P_e(x)) ,
        \left[ \proj_{Q_e} + (\proj _{\overline{\Gamma}}
    -\proj_{Q_e} ) \right] \grad P_e(x)\rangle     \\
          &  = & \langle \grad P_e(x) ,
 \grad P_e(x)\rangle
 + \langle  \grad P_e(x) ,  (\proj _{\overline{\Gamma}}
    -\proj_{Q_e} )  \grad P_e(x)\rangle \\
  &&\quad  + \ \langle \grad \pbar (x) - \grad P_e(x)) ,  \grad
P_e(x)\rangle \\
  &&\quad  + \ \langle \grad \pbar (x)
 - \grad P_e(x)) , (\proj _{\overline{\Gamma}}
    -\proj_{Q_e} )  \grad P_e(x)\rangle  \\
 &\geq & \rho^2 - \big[\rho\cdot\textstyle{\frac{1}{2}} \rho +
    \frac{45}{198}\rho\cdot \rho +
    \frac{45}{198}\rho\cdot\textstyle{\frac{1}{2}}\cdot \rho\big]
    > \frac{1}{10}\rho^2 > 0.
\end{eqnarray*}
\lqed
Following Lemma \ref{lem-innerproduct}, we require $\epsilon <
\frac{1}{2250}$ and set $r=100\epsilon.$  Then the set $\nbar$
intersects each normal slice at least once in $\calt^r$.

\bigskip
\noindent
{\bf STEP 5.}  We will show that each zero of
$\doub$ is nondegenerate with degree
1 when restricted to a normal slice.  It will follow that $\nbar$ is a
smooth section of the retraction $\pi_{N_e}$, with the $C^0$ distance from
$N_e$ to $\nbar$ less than $100\epsilon.$  The estimate used to
establish nondegeneracy will also give us the desired bound on the
$C^1$ distance $d(N,N_e)$.  For our considerations in this section, we
require $\epsilon <\frac{1}{20000}$.

At each point $x$ of $\nbar$, the projected gradient
$\doub$, being zero at that point,
has an {\bf intrinsic derivative} $\cald \doub$ which is a linear
map from $T_x M$ to the quasi-vertical space $Q_e$ at $x$.  We will show
that $\cald \doub$ is surjective with positive determinant when
restricted to $Q_e$.  From this, it will follow from the implicit
function theorem that $\nbar$ is a smooth manifold which is
the graph of a smooth section of the retraction $\pi_{N_e}$.   By our
construction in Step 5, the $C^0$ distance from
$N_e$ to $\nbar$ less than $100\epsilon.$  In fact, since any element
of $G$ could have been chosen as $e$, and the definition of $\nbar$ is
independent of this choice, the $C^0$ distance from each
$N_g$ to $\nbar$ less than $100\epsilon.$

Since the
tangent space to $\nbar$ at $x$ is the kernel of $\cald \doub$,
we will be able to estimate the $C^1$ distance by comparing the sizes
of the restrictions of $\cald \doub$ to $Q_e$ and its orthogonal complement.
Although $\cald \doub$ is defined only at the zeros of $\doub$, it
agrees there with an object defined throughout $\calt^r$, namely the
covariant derivative $\covder \doub$, considered as the operator
$v\mapsto \covder_v \doub$.  (Although
the values of $\covder \doub$ do not generally lie in $Q_e$, they do
lie in there at the zeros of $\doub$, as we shall see.)

To simplify
notation, we will denote the operation of averaging over $G$ by
$\ave$.  Then
\begin{eqnarray*}
\label{eq-covdoub}
\covder \doub & = &\covder  \proj_{Q_e}\proj_{\overline{\Gamma}}
\grad \ave P_g \\
& = & (\covder \proj_{Q_e}) \proj_{\gammabar}\ave \grad P_g
 + \  \proj_{Q_e} (\covder \proj_{\gammabar}) \ave\grad P_g
+  \proj_{Q_e} \proj _{\gammabar} \covder \ave\grad P_g  \ .
\end{eqnarray*}

We analyze the last expression above term by term.

The first term,
$ (\covder \proj_{Q_e}) \proj_{\gammabar} \grad \ave P_g,$
 is zero along the zero set $\nbar$ of the projected gradient
$\proj_{\gammabar} \grad \ave P_g=\proj_{\gammabar} \grad \pbar.$

For the second term, we have $\|\proj_{Q_e} (\covder
\proj_{\gammabar}) \ave \grad P_g \|
\leq \|\covder \proj_{\gammabar}\| \sup_g \|\grad P_g\|.$  By
Proposition \ref{prop-covderproj}, this is bounded above by $11$ times
the distance from $N_g$, i.e. by $11 \cdot 100\epsilon =1100\epsilon.$

In the last term, $ \proj_{Q_e} \proj _{\gammabar} \covder \ave \grad
P_g =\ave \proj_{Q_e} \proj _{\gammabar} \covder \grad P_g,$ the
factor $\covder \grad P_g$ is the hessian $H_{N_g}$, which we will
denote by $H_g$.  Using Proposition \ref{prop-hestimate}, we will
estimate the expression $ \proj_{Q_e} \proj _{\gammabar} \covder_v
\grad P_g$ by looking at the inner product
$\langle  \proj_{Q_e} \proj _{\gammabar} \covder_v
\grad P_g,w\rangle = H_g(v,\proj_{\gammabar}w),$ where $w$ is an
arbitrary quasi-vertical vector.   We analyze separately the cases
where $v$ is quasi-vertical and quasi-horizontal.  Since Proposition
\ref{prop-hestimate} uses the splitting into vertical and horizontal
spaces, we have to estimate several ``correction terms.''

\begin{lemma}
\label{lem-qv}
For any quasi-vertical vector $v$ at any point of $\nbar$,\newline
$\langle \cald \calv (v),v\rangle \geq \frac{89}{200} \|v\|^2.$
\end{lemma}
\pf
For $v$ quasi-vertical,
\begin{eqnarray}
\label{eq-fourterms}
  H_g(v,\proj_{\gammabar} w) &=&
   H_g(\proj_{Q_e} v,\proj_{\gammabar} w) \\  \nonumber
& =&
    H_g(\proj_{\Gamma_g} v, \proj_{\Gamma_g}w)
+   H_g((\proj_{Q_e} - \proj_{\Gamma_g})v, \proj_{\Gamma_g}w) \\ \nonumber
&&\quad +   H_g(\proj_{\Gamma_g}v,(\proj_{\gammabar} - \proj_{\Gamma_g})w)
+  H_g((\proj_{Q_e} - \proj_{\Gamma_g})v,(\proj_{\gammabar} -
\proj_{\Gamma_g})w).
\end{eqnarray}

We estimate first the differences between projections.  By Proposition
\ref{prop-twobundles} and Proposition \ref{prop-verticalangle},
\begin{eqnarray*}
\|\proj_{Q_e} - \proj_{\Gamma_g} \| &\leq &
  \|  \proj_{Q_e} - \proj_{\Gamma_e} \| + \|\proj_{\Gamma_e} -
\proj_{\Gamma_g}\|
  <   \textstyle{\frac{1}{5}} d(x,N_e)^2 + 250\epsilon \\
& < & \textstyle{ \frac{1}{5}} \big(\frac
{100}{20000}\big)^2 +  \ \frac{250}{20000} < \frac{1}{50} \ .
\end{eqnarray*}
  By Lemma
\ref{lem-closeforgauss}, $\|\proj_{\gammabar} - \proj_{\Gamma_g}\| <
1000\epsilon < \frac{1000}{20000} =\frac{1}{20}.$
  From this we have in particular that
\[
\|\proj_{\Gamma_g} v\| \geq
(1-\|\proj_{\Gamma_g} -
\proj_{\Gamma_e}\|)~\|v\| > \left(1-\frac{250}{20000}\right)\|v\| >
\frac{49}{50}
\|v\| \ .
\]

Applying Proposition \ref{prop-hestimate} now gives the estimate
$H_g(\proj_{\Gamma_g} v,\proj_{\Gamma_g} v) > .64 (\frac {49}{50})^2
\|v\|^2 > \frac{3}{5} \|v\|^2.$   Since, by \ref{prop-hestimate},
$\|H_g\|<1.32,$
the remaining three  terms in the
last expression in Equation \ref{eq-fourterms} are bounded above in
norm by $1.32\|v\|^2$ times $\frac{1}{50},$ $\frac{1}{20},$ and
$\frac{1}{50}\frac{1}{20}$ respectively, so that
$H_g(v,\proj_{\gammabar} v) \geq \frac{1}{2} \|v\|^2$ for each
quasi-vertical vector $v$.

We return now to Equation \ref{eq-covdoub} to estimate the derivative
of $\calv$ in the quasi-vertical direction.    From that equation, we
know that, for quasi-vertical $v$,
$\langle \cald \calv (v),v\rangle = \langle \covder _v \calv \rangle$
is the sum of three terms, the first of which vanishes and second of
which is bounded in absolute value by $1100\epsilon \|v\|^2$ along
$\nbar$.  The third term is the average over $G$ of quantities which
we have just shown to be bounded below by $\frac{1}{2} \|v\|^2.$
Since $\epsilon < \frac{1}{20000},$ we have
$\langle \cald \calv (v),v\rangle \geq \frac{89}{200} \|v\|^2$ for all
quasi-vertical $v$.
\lqed
  It follows from Lemma \ref{lem-qv} that the derivative $\cald \calv$
is invertible on the quasi-vertical space, so that all the zeros of
$\calv$ are nondegenerate along the normal slices.  Moreover
the positivity of $\langle \cald \calv (v),v\rangle$ shows that the
degree of each zero is $1$, so there is exactly one zero along each
such slice.
Thus, $\nbar$ is the image of a
smooth section of the normal bundle to $N_e$, i.e. of the tubular
neighborhood retraction $\pi_{N_e}$.  To estimate the derivative of
this section, and hence the $C^1$ distance from $N_e$ to $\nbar$, we
need to estimate $\cald \calv$ when applied to quasi-horizontal vectors.

\begin{lemma}
\label{lem-qh}
For any quasi-horizontal vector $v$ at any point
of the zero set $\nbar$ of $\calv$, $\|\cald \calv (v)\| \leq
60\sqrt\epsilon~ \|v\|$.
\end{lemma}
\pf We will estimate the inner product $\langle \cald\calv (v), w \rangle
= $ for quasi-horizontal $v$ and quasi-vertical $w$.  Since we already
have the bound of $1100 \epsilon \|v\|~\|w\|$ for what corresponds to
the first two terms of Equation \ref{eq-covdoub}, what remains to be
estimated is $H_g(v,\proj_{\gammabar}w)$.  For this we use Corollary
\ref{cor-hv}, making the usual slight adjustment between horizontal
and quasi-horizontal vectors, analogous to Equation
\ref{eq-fourterms}.
\begin{eqnarray}
\label{eq-fourtermsbis}
  H_g(v,\proj_{\gammabar} w) & = & H_g(\proj_{Q^{\perp}_e}
   v,\proj_{\gammabar} w) \\ \nonumber & =& H_g(\proj_{\Gamma^{\perp}_g}
   v, \proj_{\Gamma_g}w) + H_g((\proj_{Q^{\perp}_e} -
   \proj_{\Gamma^{\perp}_g})v, \proj_{\Gamma_g}w) \\ \nonumber
&&\qquad +
   H_g(\proj_{\Gamma^{\perp}_g}v,(\proj_{\gammabar} -
   \proj_{\Gamma_g})w) + H_g((\proj_{Q^{\perp}_e} -
   \proj_{\Gamma^{\perp}_g})v,(\proj_{\gammabar} -
   \proj_{\Gamma_g})w).
\end{eqnarray}
We estimate this expression term by term, as in the proof of Lemma
\ref{lem-qv}, using Corollary \ref{cor-hv} (and the duality of
Corollary \ref{cor-duality}):
$|H_g(v,\proj_{\gammabar} w) | \leq [3\sqrt{r} +1.32 (\frac{1}{5}r^2
+250\epsilon +  1000\epsilon + ( \frac{1}{5}r^2
+250\epsilon)1000\epsilon)]~ \|v\|~\|w\|.$  Using the assumptions
$r=100\epsilon$ and $\epsilon < \frac{1}{20000}$ and some arithmetic
we get $H_g(v, \proj_{\gammabar} w) \leq (30\sqrt{\epsilon}
+1670\epsilon)\|v\|~\|w\|.$

As with the previous lemma, we return to Equation \ref{eq-covdoub}
which expresses $\cald \calv (v) (v) = \covder _v (v)$ as a sum
of three terms, the first of which vanishes and second of which is
bounded in norm by $1100\epsilon \|v\|$ along $\nbar$.  The
third term is the average over $G$ of operators which we have just
shown to be bounded above in norm by $(30\sqrt{\epsilon}
+1670\epsilon)\|v\|$.  As a result, we find
$\|\cald \calv (v)\| \leq (30\sqrt{\epsilon} +2770\epsilon)\|v\|\leq
60 \sqrt {\epsilon}  \|v\|$.
\lqed

Combining the results of Lemmas \ref{lem-qv} and \ref{lem-qh}, we
conclude that the operator from $Q^{\perp}_e$ to $Q_e$
whose graph is the tangent space $T_x\nbar$ has norm
at most $60\frac{200}{89} \sqrt\epsilon < 135\sqrt{\epsilon}$.
By Corollary \ref{cor-graph}, $d(T_x \nbar,Q^\perp_e) < \tan^{-1}
( 135\sqrt{\epsilon}) <  135\sqrt{\epsilon}$.  Combining this estimate
with Proposition \ref{prop-twobundles}, we find that the Finsler distance from
$\nbar$ to the parallel translate $\Gamma^{\perp}_e$ of the tangent
space to $N_e$ at the foot of the normal from $x$ is less than
$136\sqrt{\epsilon}.$  This Finsler distance is one of the two quantities
entering
in the definition of the $C^1$ distance from $N_e$ to $\nbar$.  The
other is simply $d(x,N_e)$, which by
Step 4 is bounded above by $r=100\epsilon < \sqrt{\epsilon}$.
Since $e$ could have been any element of $G$, we conclude that
$d(N_g,\nbar) < 136 \sqrt{ \epsilon}$ for all $g$ in $G$, and
Theorem \ref{thm-centerofmass} is proven.
\qed

\section{Other ways of averaging submanifolds}
\label{sec-other}
Before (and while) obtaining the present proof of Theorem
\ref{thm-centerofmass}, we tried several other ways to average
submanifolds.  None of them have succeeded yet, but there may be means
of circumventing the difficulties which we found.  In each subsection
below, we will
describe a proposed procedure for averaging a pair of
submanifolds; we leave it to the reader to general each one to
arbitrary weighted families.

\subsection{Pursuit}
Given submanifolds $N_1$ and $N_2$, we may move each point on $N_i$
half-way along the geodesic segment
to the nearest point on the other submanifold.  The resulting
two submanifolds should be closer together than the original pair, and
we may then iterate the process and hope that the pair converges to a
single submanifold.
Variants on this method would involve going only a small fraction of
the way from each submanifold to the other at each step, or letting each
submanifold
evolve continuously by the (time dependent) vector field pointing toward the
nearest points on the other one.

We did not succeed in obtaining the a priori bound on derivatives
which would insure convergence for any version of this procedure.  It
might be interesting to test the idea with numerical experiments.

\subsection{The space of submanifolds}
One could consider $N_1$ and $N_2$  as points in the space of
submanifolds, and then apply Grove-Karcher averaging there.  (For a
pair of submanifolds, this would just mean going half-way along the
geodesic from one manifold to the other.)  The problem here is that
there is no natural differentiable structure on a Banach manifold of
unparametrized submanifolds.  In fact, if one uses the natural
coordinate systems given by sections of normal bundles, the transition
maps fail to be differentiable at points which do not have extra
smoothness.
This difficulty is related to the fact that
multiplication and inversion are not differentiable maps for any
Banach manifold of diffeomorphisms.  Perhaps there is still a
tractable variational problem for paths in the space of submanifolds
for which two nearby points are connected by a unique minimizer.
 (See Section \ref{subsec-transport}
for a related remark.)

\subsection{Distribution theory}
The average of $N_1$ and $N_2$ is well defined as a de Rham current on
$M$, i.e. as a continuous linear functional on the differential forms
of degree equal to the dimension of the  $N_i$ (at least if the submanifolds
are orientable).  This current is closed and is ``approximately a
manifold'' in some sense.  Perhaps there is a theorem which allows one
to ``condense'' such a current to an honest manifold.  We don't even
know how to begin proving such a theorem.

Alternatively, since the submanifolds inherit riemannian metrics from
$M$, one could consider their average as a measure on $N$ and, once
again, try to turn this measure into one supported on a submanifold.
This idea seems even harder to implement than the one using currents,
since one does not easily read the differentiable structure of the
submanifolds from the measures.

\section{Final remarks}
\label{sec-final}

\subsection{Is the average gentle?}

Unfortunately, our method gives no information on the extrinsic
curvature of $\nbar$, so cannot establish that $(M,\nbar)$ is a gentle
pair.  In particular, we have no way to estimate the opposite
distances $d(\nbar,N_g)$.  It would be interesting to know whether the
averaging of gentle pairs (even gentle curves in euclidean 3-space)
could produce a highly ``wrinkled'' submanifold.  Perhaps one of the
alternative averaging methods suggested in Section \ref{sec-other}
would be preferable in this respect.  

We note that, for a pair of hypersurfaces with equal weights, our
averaging method just produces the hypersurface of points which are
equidistant from the original two.  In this special case, there is no
``loss of regularity,'' because the average hypersurface can be found
directly from the potential functions of the original two, rather than
from their gradients.  

\subsection{Possible extensions of the theorem}

Is there a more global version of Theorem \ref{thm-main}, which
would subsume Cartan's fixed point theorem (see the Introduction).
In addition to the
assumption of nonpositive curvature, one would presumably need to impose
conditions on the submanifold $N$, for instance that its normal
exponential map be a diffeomorphism.

It should be possible to extend our theorem to the case where the
ambient manifold $M$ is infinite-dimensional.  This extension would be
essential for applications to almost-actions of non-finite compact Lie
groups.  (See the Introduction.)  If $M$ is a Hilbert
manifold, the only place where more work needs to be done is in Steps
4 and 5 of the proof of the averaging theorem, where degree theory is
used.  (See Remark \ref{rmk-infinite}.)  If $M$ is a general Banach
manifold, we would also need to extend our study of the Grassmann
manifolds to this case.  This would be complicated by the fact that
orthogonal projections are no longer available.

There should also be versions of our theorem in symplectic geometry,
e.g. for isotropic, lagrangian, or coisotropic submanifolds.   If
there are not, the obstructions to solutions should be interesting.

\subsection{Morphing}
Jamie Sethian pointed out that a solution to our averaging problem for
a pair $N_1$ and $N_2$ of submanifolds would, on repeated application,
give a family $N_t$ of submanifolds defined for all binary fractions
$t \in [1,2]$.  Presumably, this could be extended by continuity to
all real $t \in
[1,2]$, thus giving a deformation, or ``morphing'', between the two
submanifolds.  (Alternatively, we could vary the weight on the
parameter space $\{1,2\}$ to move from one submanifold to the other.)
Conversely, any ``reversible'' procedure for morphing one
manifold
into another would give a solution to the averaging problem.

\subsection{Mass transport}
\label{subsec-transport}
If the two submanifolds $N_i$ are parametrized by the same manifold $N$,
the midpoint problem
becomes trivial, since we can simply take the midpoints of the
geodesic segments joining correspondingly parametrized points.  (Of
course, some estimates would still be needed to show that the
resulting set of midpoints is a manifold.)  Conversely, a solution of
the midpoint problem gives a diffeomorphism between $N_1$ and $N_2$,
obtained by projecting the average manifold onto both of them. 
David Aldous pointed out to me that this problem of finding a ``distinguished''
diffeomorphism between two submanifolds may be thought of as a special
case of the Monge-Kantorovich mass-transport problem (see, for
example, \cite{ga-mc:geometry}), in which an optimal correspondence is
sought between the supports of two measures on the same space.
Unfortunately, the correspondences found in most known solutions of
this problem are multiple-valued.

\subsection{Evolution of submanifolds}
\label{subsec-evolution}
The representation of submanifolds as zeros of projected gradient
functions may be useful for describing the evolution of manifolds under
flows, generalizing the level set method for hypersurfaces.  Ruuth et
al \cite{ru-me-xi-os:diffusion} have studied such evolutions
numerically.  On the other hand, Ambrosio and Soner \cite{am-so:level}
study these evolutions by representing a submanifold as the zero set
of its distance or potential function.  Although they use a single
function and not a section of a bundle, they must analyze the
derivatives of the potential function and the curvature of tube
boundaries just as we do.

\appendix
\section
{Appendix: Finsler distance and averaging in Grassmann manifolds}
\label{sec-appendix}
We review in this section some geometry of Grassmann manifolds.
Specifically, we will describe a center-of-mass construction
derived from the averaging of projections, and we will analyze this
construction in terms of a natural Finsler metric on the grassmannian.
Many of the results here are probably not new, but there seems to be
no convenient reference for them in the form which we need.

We will work with finite dimensional subspaces of (possibly infinite
dimensional) Hilbert spaces.  By duality, we will extend our results
to subspaces of finite codimension (such as the normal spaces to
finite dimensional submanifolds), but this is as far as we can go.
The study of grassmannians of general subspaces is much more subtle.
(See, for example, Chapter 7 of \cite{pr-se:loop}.)

For our purposes, a Finsler manifold will be a Banach manifold with a
continuous field of norms on the tangent spaces.  The norms should
define the topology but need not be smooth or strictly convex.  These
conditions are enough to define length of curves and an inner metric
(see \cite{ri:innere}, \S16)
compatible with the topology, though nearby points may
be connected by many length-minimizing paths, such as in the case of
the $L^{\infty}$ (``taxicab'') metric on the plane.

\subsection{The Grassmannian and its metric}

For a real Hilbert space $E$ and $k=0,1,2,\ldots$
we denote by $G_k(E)$ the Grassmann
manifold of $k$-dimensional subspaces of $E$, and by $G(E)$ the
manifold having the $G_k(E)$'s as its connected components.
Each $G_k(E)$ is a homogeneous space under the
natural action of the orthogonal group $O(E)$; the isotropy of an
element $F\subset E$ is the direct product $O(F) \times O(F^\perp)$.
Each tangent space $T_F G(E)$ is naturally isomorphic to the space
$\Hom (F,F^\perp)$ of continuous linear maps from $F$ to $F^\perp$.

Endowing each tangent space $\Hom(F,F^\perp)$ with the operator norm
gives an $O(E)$-invariant Finsler metric on $G(E)$.  We will use
this metric rather than the riemannian metric (also $O(E)$-invariant)
coming from the inner products $\langle Q_1,Q_2 \rangle = \Tr(Q_1
Q_2^*)$ on the tangent spaces.  The Finsler distance $d(F,F')$
between two subspaces is
then defined as the infimum of the lengths of paths joining $F$ to
$F'$.

We embed $G(E)$ into the space $S(E)$ of self-adjoint endomorphisms of
$E$ by mapping each $F$ to the operator $P(F)=\proj_F$ of orthogonal projection
onto $F$.  The image of $G_k(E)$ under this embedding consists of
the self-adjoint projections with trace $k$.  The tangent space to
$P(G(E))$ at $P(F)$ consists of those $q$ in $S(E)$ such that $\proj_F
q = q(I-\proj_F)$.  The image of $Q \in T_F G(E)$ under the derivative
$T_F P$ has the block
diagonal form
$$\displaystyle
q = \left[
\begin{array}{cc}
 0 & Q^*\\
 Q & 0
\end{array}
\right]\
$$
with respect to the decomposition $E=F\oplus F^\perp$.

The
Hilbert norm $\sqrt{\Tr (q^* q)}$ of $q$ is $\sqrt{2}$ times that of $Q$, but
the operator norms of $q$ and $Q$ are exactly the same.  Thus $P$ is
an isometric embedding of the Finsler manifold $G(E)$ into the Banach
space $S(E)$.\footnote{The map $P$ keeps the length of curves fixed,
  but it generally decreases distances; see Proposition
  \ref{prop-finslervsnorm} below.}

We will also make use of the antisymmetric operator
$$\displaystyle
h(Q) = \left[
\begin{array}{cc}
 0 & -Q^*\\
 Q & 0
\end{array}
\right]\
$$
in the Lie algebra $\mathfrak{o}(E)$ of $O(E)$.  We call it the {\bf
  horizontal lift} of the tangent vector $Q$, because the horizontal
lifts of all $Q \in T_F(G(E))$ form a complement to the isotropy
subalgebra $\mathfrak{o}(F) \times \mathfrak{o}(F^\perp)$ of $F$.  By
right translating the horizontal space around $O(E)$, we obtain a
connection in the principal bundle $O(E)\to G_k(E)$.    Thus, each
path $\{F_t\}$ in $G(E)$ has a horizontal lift $\{g_t\}$ in $O(E)$
such that $F_t=g_t(F_0)$, and $dF_t/dt$ and $dg_t/dt$ have the same
norm.

It will be useful to compare the Finsler distance with some other
notions of distance between subspaces.  We have the following results.

\begin{prop}
\label{prop-finslervsC0}
If $F$ and $F'$ are subspaces of the same (finite) dimension in
$E$, then the Finsler distance $d(F,F')$
is equal to the
$C^0$ distance $d_0(F,F')$
between the unit spheres of $F$ and $F'$, considered
as submanifolds of the unit sphere in $E$ (with the standard riemannian
metric).
\end{prop}

\begin{cor}
The Finsler distance between any two elements of $G_k(E)$ is at most
$\pi/2$.
\end{cor}

\begin{cor}
The $C^0$ distance $d_0(F,F')$ is symmetric in $F$ and $F'$.
\end{cor}

\begin{prop}
\label{prop-finslervsnorm}
If $F$ and $F'$ are subspaces of the same (finite) dimension in
$E$, then $\| \proj_{F} - \proj_{F'}\| =\sin d(F,F').$
\end{prop}

\bigskip
\noindent
{\bf Proof of Proposition \ref{prop-finslervsC0}.}
To begin, we note that, by an arbitrarily small perturbation of $F$,
we can insure that
$$F \cap F'^\perp = \{ 0\}.~~~~~~~~~~~~~~~~~~(*) $$  It will suffice,
then, to prove our equality for pairs having property $(*)$.  (When it
is not satisfied, $d_0(F,F')=\pi/2$.)

 We will show first that $d(F,F') \leq d_0(F,F')$.
By our assumption $(*)$, with respect to the decomposition $E=F\oplus
F^\perp$, $F'$ is the graph of a linear map $\phi:F\to F^\perp$.
The non-negative symmetric operator $\phi^* \phi$
on $F$ has $k$ eigenvalues which we take to be the squares of
 $$\lambda_1 \geq
\ldots \geq \lambda_k \geq 0.$$
We denote an orthonormal basis of
eigenvectors with these eigenvalues by $e_1,\ldots,e_k$, and we call
the angles $\theta_j = \tan^{-1} \lambda_j$ the {\bf canonical
  angles}\footnote{These angles  have a long history, going
  back to Jordan \cite{jo:essai}.  We refer to \cite{pa-we:history}
  for an extensive discussion of this history, and add as well the
  references \cite{af:orthogonal} and \cite{pe:integralgeometrie}.
  Finally, we note that Porta and Recht \cite{po-re:minimality} have
  quite explicitly studied some aspects of the Finsler geometry of
  grassmannians.  They represent subspaces by the involutions
  $2\proj_F - I$, so their Finsler metric is twice ours.}
between the subspaces $F$ and $F$'.

Since
$$\langle \phi e_j,\phi e_k \rangle = \langle \phi^*\phi e_j,e_k
\rangle = \lambda_j^2 \delta _{jk}, $$
we get an orthonormal basis for $F'$ by taking the vectors
$$\frac{e_j + \phi e_j }{\sqrt{1+\lambda_j ^2}} = \frac{e_j + \phi e_j
}{\sqrt{1+\tan^2 \theta_j}}=( \cos\theta_j) (e_j + \phi e_j).$$
We get an orthonormal set in $F^\perp$ by taking $f_j = (1/\lambda_j) \phi
e_j$ whenever $\lambda_j \neq 0$.  For the remaining values of $j$, we
set $f_j = 0$.

Our basis elements  for $F'$
can now be written as $\cos \theta _j e_j + \sin \theta _j f_j$, and
we can now  connect $F$ and $F'$ by the path
$$F_t= \cos (t\theta_j)e_j +\sin ( t\theta_j)f_j, ~~~~0 \leq t \leq 1.$$
The velocity of this path at $t=0$ is the map $e_j\mapsto \theta_j
f_j$, which has norm $\theta_1$, the largest canonical angle.  Since
the $F_t$ form a trajectory starting at $F_1$ for a 1-parameter group
of isometries (generated by the skew-symmetric transformation which
takes $e_j$ to $\theta_j f_j$ and is zero on $(F \oplus \phi F)^\perp$),
the speed of the path is constant, so its length is
$\theta_1$.  Hence $d(F,F')\leq \theta_ 1$.
But $\theta_1 \leq d_0(F,F') $ because $\theta _1$ is the distance from
$\cos \theta _1 e_1 + \sin \theta _1 f_1 \in F'$ to the unit sphere in $F$, so
we have proven that $d(F,F')\leq d_0(F,F')$.

For the reverse inequality, we will show that $d_0(F,F')$ is bounded
above by the length of any path from $F$ to $F'$.  Let $\{F_t\}$ be
such a path, $\{g_t\}$ its horizontal lift to the orthogonal group.
Any unit vector in $F'$ is $g_1 v$ for some unit vector $v$ in $F$; the
distance from $g_1 v$ to the unit sphere in $F$ is estimated by using the path
$\{g_t v\}$:
$$d(v,g_1 v) \leq \int_0^1 \|\frac{d}{dt} (g_t v)\|~dt \leq
\int_0^1 \|\frac{d}{dt} g_t\|~dt = \int_0^1 \|\frac{d}{dt} F_t~dt\|,$$
which is the length of the path  $\{F_t\}$.
\qed

\bigskip
\noindent
{\bf Proof of Proposition \ref{prop-finslervsnorm}.}
Let $\theta_1 \geq \ldots \geq \theta_k$ be the canonical angles between
$F$ and $F'$.  According to the preceding proposition,
$d(F,F')=\theta_1$.  In an orthonormal basis beginning with
the orthonormal set $e_1,f_1,e_2,f_2,\ldots$
constructed in the proof above, the projection
 $\proj_{F'}$ is block
diagonal with $2\times 2$ blocks of the form
$$\displaystyle
\left[
\begin{array}{cc}
\cos^2 \theta_j & \cos \theta_j \sin \theta_j\\
 \cos \theta_j \sin \theta_j & \sin^2 \theta_j
\end{array}
\right]\
$$
followed by zero blocks,
while $\proj_F$ has the same form with  the
$\theta_j$'s replaced by 0.  The difference $\proj_F - \proj_{F'}$ is
built of diagonal blocks
$$\displaystyle
h(Q) = \left[
\begin{array}{cc}
1-\cos^2 \theta_j &- \cos \theta_j \sin \theta_j\\
 -\cos \theta_j \sin \theta_j & -\sin^2 \theta_j
\end{array}
\right] = \sin \theta_j \left[
\begin{array}{cc}
\sin \theta_j &- \cos \theta_j\\
 -\cos \theta_j  & -\sin \theta_j
\end{array}
\right],\
$$
each of which has norm $\sin \theta_j$, so
$$\|\proj_F-\proj_{F'}\| = \sin \theta_1 = \sin d(F,F').$$
\qed

The identification of $d(F,F')$ with the largest canonical angle also
leads immediately to the following result.
\begin{cor}
\label{cor-graph}
If $F'$ is the graph of a linear map $u$ from $F$ to $F^{\perp}$, then
$d(F,F')=\tan^{-1} \|u\|.$
\end{cor}

It is easy to go back and forth between subspaces and their
orthogonal complements.
\begin{cor}
\label{cor-duality}
If $F$ and $F'$ are elements of $G_k(E)$, and $E$ is finite
dimensional, then $d(F,F') = d(F^{\perp},F'^{\perp}).$
\end{cor}
\pf Note that the projections on a subspace and its orthogonal complement
sum to the
identity operator.  Now apply Proposition \ref{prop-finslervsnorm}.
\qed

\begin{rmk}
{\em
Thanks to Corollary \ref{cor-duality} (or simply using Proposition
\ref{prop-finslervsnorm} as a definition of the distance), we can
define the distance between subspaces of finite {\em codimension} in a
Hilbert space to be the distance between their orthogonal complements.

We also note that our distance functions are
unaffected by rescaling of the inner product on $E$.
}
\end{rmk}

\subsection{Center of mass in the
Grassmannian}

We turn now to the averaging of subspaces.  Using Grove-Karcher
averaging in the grassmannian makes it more difficult to
obtain estimates of the kind we
need, since it is hard to control the behavior of the Finsler distance
along geodesics.
Instead, we average the orthogonal projections associated to the
subspaces.  This approach was used already by de la Harpe and Karoubi
in \cite{de-ka:representations}, but their work is more
complicated and their estimates less sharp than ours because their
projections are not necessarily orthogonal.

Suppose, then, that we are given a family $F_g$ of elements of
$G_k(E)$ parametrized measurably by elements $g$ of a probability
space $G$.  If we average the projections $\proj_{F_g}$, we obtain a
self-adjoint operator which is not, in general, a
projection.  However, if the $F_g$ are sufficiently close together,
this operator is an ``almost projection,'' and we will see that it can be
approximated in a canonical way by a self-adjoint projection having
rank $k$.   Here is the precise result.

\begin{prop}
\label{prop-averageproj}
Let $F_g$ be a family  of elements of
$G_k(E)$ parametrized measurably by elements $g$ of a probability
space $G$ such that $d(F_g,F_h) < \epsilon$ for all $g$ and $h$ in
$G$.   As long as $\epsilon < \frac{1}{2}$, we can define an ``average''
element $\overline{F} \in G_k(E)$ such that $d(F_g,\overline{F})
< \sin^{-1} 2\epsilon$ for all $g$ in $G$.
This averaging construction is equivariant with respect to
isometries of $E$ and measure preserving automorphisms of $G$.  If the
$F_g$ also depend smoothly on other parameters $\mu$, the average
$\overline{F}$ also depends smoothly on these parameters.
Furthermore, if $\epsilon < \frac{1}{4}$, the magnitude of the
derivative $\frac{\partial \fbar}{\partial \mu}$ is at most 8 times
the supremum over $g$ of the magnitudes of the
 derivatives $\frac{\partial F_g}{\partial \mu}$.
\end{prop}
\pf
We note first of all that, by Proposition \ref{prop-finslervsnorm},
$\|\proj_{F_g} - \proj_{F_h}\| < \epsilon.$
For the average
$$\overline{P}=\int_G \proj_{F_g}~dg$$ in the space of
self-adjoint operators on $E$, a simple computation shows that
$\|\overline{P} - \proj_{F_g}\| < \epsilon$ for all $g$ in $G$.

We will now modify the approximate projection $\overline{P}$ to make
it into an orthogonal projection.  The estimate
$\|\overline{P} - \proj_{F_g}\| < \epsilon$
 tells
us that the spectrum of $\overline{P}$ is contained in intervals of radius
$\epsilon$ about $0$ and $1$.  Since $\epsilon < \frac{1}{2},$ we can
apply to $\pbar$ the holomorphic function $h$, defined on the complex
$\lambda$ plane with the line $\rm{Re}~ \lambda = \frac{1}{2}$ removed,
having the value $0$ to the left of the line and $1$ to the right.
The operator $h(\pbar)$ is a spectral projection of $\pbar$ which
can be written as an integral around a circle $\gamma$ of radius
$\frac{1}{2}- \delta$ around $\lambda = 1$, where $\delta$ can be an
arbitrarily small positive number.
$$ h(\pbar) = \frac{1}{2\pi i} \int_{\gamma} (\lambda I -
\pbar)^{-1}~d\lambda.$$  The application of $h$ to $\pbar$ amounts to
moving the spectrum to $0$ and $1$ while keeping the invariant
subspaces the same, so that $\|h(\pbar)-\pbar\| < \epsilon.$

Denoting by $\fbar$ the image of $h(\pbar)$, we have $h(\pbar)$ =
$\proj_{\fbar}$, and $\|\proj_{\fbar}- \proj_{F_g}\| \leq
\|h(\pbar)-\pbar\| + \|\pbar - \proj_{F_g}\| < \epsilon + \epsilon =
2\epsilon < 1$ for all $g$ in $G$.  Since $\|\proj_{\fbar}-
\proj_{F_g}\| < 1$, the spaces $\fbar$ and $F_g$ have the same
dimension (see \cite{ka:perturbation}, Theorem 6.32), so $\fbar \in
G_k(E)$.  By
Proposition \ref{prop-finslervsnorm}, $d(\fbar,F_g) < \sin^{-1} 2\epsilon$.

The equivariance of our construction is obvious, since all the steps
were canonical.  To establish smooth dependence and estimate
derivatives with respect to parameters, we continue our analysis of
the operator function $h$ by observing first that differentiation
under the integral sign gives smoothness (in fact analyticity) of $h$.
To estimate the derivative, we assume that the spectrum of an operator
family $P(\mu)$ lies in the union of discs of radius $\epsilon$ around
$0$ and $1$ and compute:
\begin{eqnarray*}
\frac{\partial h(P(\mu))}{\partial \mu}& =& \frac{1}{2 \pi i}
\int_{\gamma} \frac{\partial}{\partial \mu} (\lambda I - \pbar)^{-1}
~d\lambda \\
& = & - \frac{1}{2 \pi i}\int_{\gamma} (\lambda I - P(\mu))^{-1}
\frac{\partial P(\mu)}{\partial\mu} (\lambda I - P(\mu))^{-1}~d\lambda.
\end{eqnarray*}
Hence,
$$ \|\frac{\partial h(P(\mu))}{\partial \mu} \| \leq
\frac{1}{2\pi} (\rm{length}~ \gamma) \sup_{\lambda} \|(\lambda I -
P(\mu))^{-1}\| ^2 ~\|\frac{\partial P(\mu)}{\partial \mu}\|.$$
As the radius of $\gamma$ approaches $\frac{1}{2}$,
the length of $\gamma$ approaches $\pi$, and so
$$\|\frac{\partial h(P(\mu))}{\partial \mu} \| \leq\frac{1}{2\pi} \pi
\frac{1}{(\frac{1}{2}-\epsilon)^2}  \|\frac{\partial P(\mu)}{\partial
\mu}\|
\leq \frac{2}{(1-2\epsilon)^2} \|\frac{\partial P(\mu)}{\partial\mu}\|.
$$
As long as $\epsilon < \frac{1}{4},$ we have $$\| \frac{\partial
h(P(\mu))}{\partial \mu} \| \leq 8  \|\frac{\partial P(\mu)}{\partial
\mu}\|,$$ and we shall be content with this crude estimate (though the
factor $8$ can probably be brought close to unity as $\epsilon$ becomes
small).

We can now finish the proof of our Proposition.  Let the spaces $F_g$
depend smoothly on $\mu$.  Since the map $F\mapsto \proj_F$ is a
Finsler isometry, $\|\frac{\partial F_g (\mu)}{\partial \mu} \|
= \|\frac{\partial \proj_{F_g}}{\partial \mu}\|.$  Averaging does not
increase norms, so $\|\frac{\partial \pbar (\mu)}{\partial \mu}\| \leq \sup_g
\|\frac{\partial F_g}{\partial \mu}\|.$  By our previous estimate,
$\|\frac{\partial \fbar}{\partial \mu}\|= \|\frac{\partial
\proj_{\fbar}}{\partial \mu} \| \leq 8~ \sup_g \|\frac{\partial
F_g}{\partial \mu}\|.$
\qed

\end{document}